\documentclass[11pt]{article}
\usepackage{float}
\usepackage{url}
\usepackage{mathrsfs}
\usepackage{color}
\usepackage[toc,page]{appendix} 

\usepackage{latexsym}
\usepackage{amsmath}
\usepackage{amsfonts}
\usepackage{amssymb}
\usepackage{amsthm}
\usepackage{psfrag}
\usepackage{epsfig}
\usepackage[english]{babel}
\usepackage{color}
\usepackage{soul}
\usepackage{marginnote}

\usepackage{tikz}
\usepackage{pgfplots}
\usepackage{pgfplotstable}
\usetikzlibrary{calc}

\newcommand{\slopeTriangleBelow}[6]
{
    \pgfplotsextra{

    \coordinate (A) at (axis cs:#1,#3);
    \coordinate (B) at (axis cs:#1,#4);
    \coordinate (C) at (axis cs:#2,#4);

    \draw[#6] (A)--(B) node [midway, left] {#5} ;
    \draw[#6] (B)--(C)  node [midway, above] {1};
    \draw[#6] (C)--(A) ;
    }
}

\newcommand{\slopeTriangleAbove}[6]
{
    \pgfplotsextra{ 
    \coordinate (A) at (axis cs: #1,#3);
    \coordinate (B) at (axis cs: #2,#3);
    \coordinate (C) at (axis cs: #2,#4);

    \draw[#6] (A)--(B) node [midway, below]{1};
    \draw[#6] (B)--(C) node [midway, right]{#5};
    \draw[#6] (C)--(A) ;
    }
}

\newsavebox{\fmbox}

\newtheorem{theorem}{Theorem}[section]
\newtheorem{lemma}[theorem]{Lemma}

\newtheorem*{remark*}{Remark}
\newtheorem*{assumption*}{Assumption}

\newcommand{\pf}{{\bf Proof: }}

\newcommand{\eps}{\varepsilon}

\newcommand{\dis}{\displaystyle}

\def\Vo{\vbox{\offinterlineskip\hbox{\kern 3pt$\scriptstyle\circ$}
\kern 1pt\hbox{$V$}}}
\def\Ho{\vbox{\offinterlineskip\hbox{\kern 3pt$\scriptstyle\circ$}
\kern 1pt\hbox{$H$}}}
\def\Wo{\vbox{\offinterlineskip\hbox{\kern 3pt$\scriptstyle\circ$}
\kern 1pt\hbox{$W$}}}

\newcommand{\R}{\mathbb{R}}

\newcommand{\caL}{{\cal L}}
\newcommand{\caC}{{\cal C}}

\newcommand{\caK}{{\cal K}}

\newcommand{\caQ}{{\cal Q}}

\renewcommand{\textbf}[1]{{\bfseries\mathversion{bold}#1}}
\newcommand{\bfu}{\boldsymbol{u}}

\newcommand{\s}{{\sigma}}
\newcommand{\sn}{{\sigma_n}}

\newcommand{\vn}{v_n}

\renewcommand{\div} {{\rm{div} \,}}
\newcommand{\GC}{\Gamma_C}

\newcommand{\GCQ}{Q_C}
\newcommand{\GCQtilde}{\tilde{Q}_C}

\newcommand{\dG}{\ {\rm d} \Gamma}
\newcommand{\dO}{\ {\rm d}\Omega}
\newcommand{\intGC}{\int_{\Gamma_C}}
\newcommand{\intGN}{\int_{\Gamma_N}}
\newcommand{\intO}{\int_{\Omega}}

\newcommand{\abs}[1]{\left|{#1}\right|}
\newcommand{\norm}[1]{\left\|{#1}\right\|}

\newcommand{\scal}[2]{\left<#1,#2\right>_{W',W}}

\newcommand\restr[2]{{
  \left.\kern-\nulldelimiterspace 
  #1 
  \vphantom{\big|} 
  \right|_{#2} 
  }}
\newcommand\restrr[2]{\ensuremath{\left.#1\right|_{#2}}}
\newcommand{\cqfd}{{$\mbox{}$\hfill$\square$}}

\usepackage{letltxmacro}
\makeatletter
\let\oldr@@t\r@@t
\def\r@@t#1#2{%
\setbox0=\hbox{$\oldr@@t#1{#2\,}$}\dimen0=\ht0
\advance\dimen0-0.2\ht0
\setbox2=\hbox{\vrule height\ht0 depth -\dimen0}%
{\box0\lower0.4pt\box2}}
\LetLtxMacro{\oldsqrt}{\sqrt}
\renewcommand*{\sqrt}[2][\ ]{\oldsqrt[#1]{#2} }
\makeatother

\renewcommand{\textbf}[1]{{\bfseries\mathversion{bold}#1}}
\newcommand{\bfz}{\boldsymbol{\zeta}}
\newcommand{\bfX}{\boldsymbol{\Xi}}

\newcommand{\bfi}{\boldsymbol{i}}

\newcommand{\bfI}{\boldsymbol{I}}

\newcommand{\bflambda}{\boldsymbol{\lambda}}
\newcommand{\bfN}{\boldsymbol{N}}
\newcommand{\bfB}{\boldsymbol{B}}

\usepackage{subfigure}

\usepackage[normalem]{ ulem }
 \usepackage{soul}

\title{\textbf{ \textit{A priori} error for unilateral contact problems with augmented Lagrange multipliers and IsoGeometric Analysis }} 
\author{Mathieu Fabre\footnote{EPFL SB MATHICSE MNS (B\^at. MA), station 8, CH 1015 Lausanne (Switzerland).}$  \ ^, $\footnote{\samepage Istituto di Matematica Applicata e Tecnologie Informatiche 'E. Magenes' del CNR
via Ferrata 1, 27100, Pavia (Italy). \qquad\qquad\qquad\qquad\qquad\qquad\qquad\qquad\qquad\qquad\qquad\qquad\qquad\qquad\qquad\qquad\qquad\qquad\qquad\qquad\qquad\qquad $ $  email: mathieu.fabre@epfl.ch.}}

\date{\today}

\begin{document}

\renewcommand{\labelitemi}{\textbullet}  
\maketitle
 
\begin{abstract}
The aim of the present work is to extend the \textit{a priori} error for contact problems with an augmented Lagrangian method. We focus on unilateral contact problem without friction between an elastic body and a rigid one. We consider the pushforward of a NURBS space of degree $p$ for the displacement and the pushforward of a B-Spline space of degree $p-2$ for the Lagrange multipliers. This specific choice of space is a stable couple of spaces. An optimal \textit{a priori} error estimate inspired from the Nitsche's method theory is provided and compared to the regularity of the solution. We perform a numerical validation with two- and three-dimensions in small and large deformations with $N2/S0$ and $N3/S1$ elements.
\end{abstract}


\section*{Introduction}
\label{sec:sec0}
\addcontentsline{toc}{section}{Introduction}

The purpose of this paper is to study theoretically a Lagrange multiplier method penalized in a consistent way, the augmented Lagrangian method. Little work has been done in this area and, as far as we know, this result is the first theoretical result of an optimal a priori estimate for the contact problem for an augmented Lagrangian method using Isogeometrical Analysis. This article is based on the article by Erik Burman and his co-authors \cite{burman-16}, where it used Brouwer's fixed point theorem to show the existence and uniqueness of the approximate solution. This recent work is inspired by the work carried out on the contact of Nitsche's contact methods, where we use the coercivity and the hemi-continuity of the operator, a method that has proven its stability, its advantages and its robustness in many cases , for example fictitious domain and in dynamic cases \cite{chouly-hild-2013,chouly-hild-renard-2014,chouly-hild-renard-2015-1,chouly-hild-renard-2015-2,fabre-pousin-renard-2016,chouly-mlika-renard-2017,franz-overview}.

The methods for contact problems have been increasingly studied in recent years and remain, especially in the industry, a central point due to its intrinsic non-linearity at the edge of contact and poor conditioning \cite{al-cu1988,laursen-02,wriggers-06,kon-sch-13}.

In order to take into account this lack of robustness and to obtain an accurate method, work using the framework of isogeometric analysis  \cite{hughes05} increases. Indeed, the geometries are precisely or exactly approximate and also smoother. Moreover, in the industry, the different geometries are initially built thanks to CAD, which uses Bezier curves and their generalizations: B-Splines and NURBS. However, the isogeometric paradigm is based on the use of these basic functions in order to discretize the partial differential equations. They have many advantages, including the use of fewer degrees of freedom in order to represent the bodies and higher approximation analysis. Isogeometric analysis methods for contact problems have been introduced in \cite{lorenzis9,Temizer11,Temizer12,lorenzis12,lorenzisreview,lorenzis15}, using Lagrange methods or augmented Lagrangian methods and also see those using primal and dual elements \cite{Wohlmuth00,hueber-Wohl-05,hueber-Stadler-Wohlmuth-08,wohl12,seitz16}. 

In this paper, we theoretically and numerically extend a Lagrange multiplier method \cite{fabre-a-priori-iga-17} using the theory present in \cite{chouly-hild-2013,chouly-hild-renard-2014,franz-overview} for Nitsche's methods. We continue with the stable choice of multiplier space proposed in \cite{brivadis15} and used in \cite{fabre-a-priori-iga-17}. For the numerical point of view, an active set strategy method is used in order to help the convergence of Newton-Raphson iterations \cite{hueber-Wohl-05,hueber-Stadler-Wohlmuth-08}.

Finally, the performance of this method will be presented on tests in small and large deformations, using an internal development code, the free library of Igatools \cite{Pauletti2015}. 

In Section 1, we will introduce Signorini's problem and the various notations. Section 2 is devoted to the description of discrete spaces and their properties. In Section 3, an a priori optimal estimate will be presented. In the last section, we will illustrate cases in small and large deformations for different types of elements.

\section{Preliminaries and notations}
\label{sec:sec1}
\subsection{Unilateral contact problem}
\label{subsec:continuous_pb}

We consider that $\Omega \subset \R^d$ ($d = $2$ \textrm{ or } 3$) is a bounded regular domain which represents the reference configuration of an elastic body. Let $\Gamma$ be the boundary domain which is split into three non overlapping parts, the Dirichlet part $\Gamma_D$ with $\textrm{meas}(\Gamma_D) >0$, the Neumann one $\Gamma_N$ and $\Gamma_C$ the potential zone of contact. The elastic body $\Omega$ is submitted to volume load $f$, to surface force $\ell$ on $\Gamma_N$  and a homogeneous Dirichlet condition at $\Gamma_D$. For the next of the section, we define our normal vector $n$ as the unit normal vector of the rigid body and by $n_\Omega$ the outward normal vector on $\Gamma$.

In the following section of the article, we denote the displacement by $u$ of the domain $\Omega$,  the linearized strain tensor by $\dis \eps (u) = \frac{1}{2} (\nabla u + \nabla u^T)$ and the stress tensor by $\dis \s  = (\s_{ij})_{1 \leq i,j \leq d} $ is given by $\dis \s (u)= A \eps (u)$, where $A=(a_{ijkl})_{1 \leq i,j,k,l \leq d}$ is a fourth order symmetric tensor verifying the usual uniform ellipticity and boundedness properties.

We decompose any displacement in $\Omega$ and any density of surface force on $\Gamma$ as a normal and a tangential components, as follows: $$ u = u_n n + u_t \qquad \textrm{and} \qquad \s (u) n = \s_n (u) n + \s_t (u) .$$

We write the classical unilateral contact problem between an elastic body and a rigid one, find $u$ such that
\begin{eqnarray} \label{eq:strong}
\begin{array}{rcll}
 \div \s (u) +f \!\!\!&\!\!\!=\!\!\!&\!\!\! 0  &\qquad \textrm{in } \Omega, \\
 \s (u) \!\!\!&=&\!\!\! A \eps (u)&\qquad \textrm{in } \Omega,\\
 u \!\!\!&=&\!\!\! 0 &\qquad \textrm{on } \Gamma_D,\\
  \s (u) {n_\Omega} \!\!\!&=&\!\!\! \ell &\qquad \textrm{on } \Gamma_N,
 \end{array}
\end{eqnarray} 
and the Signorini condition without friction at $\Gamma_C$ are:
\begin{eqnarray} \label{eq:contact_cond}
\begin{array}{rl}
 u_n \geq &\!\!\! 0  \quad (i), \\
 \s_n (u) \leq &\!\!\! 0  \quad (ii), \\
  \s_n (u) u_n = &\!\!\! 0  \quad (iii), \\
 \s_t (u) = &\!\!\! 0  \quad (iv).
\end{array}
\end{eqnarray} 
We consider the following Hilbert spaces to describe the variational formulation of \eqref{eq:strong}-\eqref{eq:contact_cond}: $$ \dis V := H^1_{0,\Gamma_D}(\Omega)^d = \{ v \in H^1(\Omega)^d , \quad v = 0 \textrm{ on } \Gamma_D \}, \quad W=\{ \restrr{v_n}{\Gamma_C}, \quad v \in V \},$$ and their dual spaces $V'$, $W'$ endowed with their usual norms and we denote by $\scal{\cdot}{\cdot}$ the duality pairing between $W'$ and $W$.

We introduce the following notations, we denote by $\norm{\cdot}_{3/2+s,\Omega} $ the norm on $H^{3/2+s}(\Omega)^d$ and by $\norm{\cdot}_{s,\GC} $ the norm on $H^s(\Gamma_C)$. For all $u$ and $v$ in $V$, we set: $$\dis a(u,v) = \intO \s(u): \eps(v)\dO \quad \textrm{and} \quad L(v) = \intO f \cdot v \dO + \intGN \ell \cdot v \dG.$$

\noindent We define the classical variational inequality of \eqref{eq:strong}-\eqref{eq:contact_cond} (see \cite{lions-magenes-72}) by finding 

\noindent $ u \in K_C := \{ v \in V, \quad v_n \geq 0 \textrm{ on } \Gamma_C\}$ such as: 
\begin{eqnarray} \label{eq:var_ineq}
\dis a(u,v-u) \geq L(v-u), \qquad \forall v \in K_C,
\end{eqnarray} 
where $K_C$ is the closed convex cone of admissible displacement fields satisfying the non-interpenetration conditions.

\noindent It is well known that a Newton-Raphson's method cannot be used directly to solve this formulation \eqref{eq:var_ineq}. One method is to introduce a Lagrange method denoted by $\lambda$, which represents the surface normal force. For all $\lambda$ in $W'$, we denote $\dis b(\lambda,v) = - \scal \lambda \vn$ and $M$ is the classical convex cone of multipliers on $\Gamma_C$:
 $$ \dis M := \{ \mu \in W', \quad \scal{\mu}{\psi} \leq 0 \quad \forall \psi \in H^{1/2}(\Gamma_C), \quad \psi \geq 0 \textit{ a.e.} \textrm{ on } \Gamma_C \} .$$
  \noindent We can now rewrite the complementary conditions as follows:
\begin{eqnarray} \label{eq:contact_cond_mult}
\begin{array}{rl}
 u_n \geq &\!\!\! 0  \quad (i), \\
\lambda \leq &\!\!\! 0  \quad (ii), \\
 \lambda u_n = &\!\!\! 0  \quad (iii).
\end{array}
\end{eqnarray} 
 
 \noindent The mixed formulation \cite{ben-belgacem-renard-03} of the unilateral contact problem \eqref{eq:strong} and \eqref{eq:contact_cond_mult} consists in finding $(u,\lambda) \in V \times M$ such that:
\begin{eqnarray} \label{eq:mixed_form}
\left\{
\begin{array}{rl}
\dis a(u,v) - b(\lambda,v) = L(v),& \dis \qquad \forall v \in V,\\
\dis b(\mu-\lambda,u) \geq 0,&\dis \qquad \forall \mu \in M.
\end{array}
\right.
\end{eqnarray} 
Stampacchia's Theorem ensures that problem \eqref{eq:mixed_form} admits a unique solution.
  
\noindent The existence and uniqueness of the solution $(u,\lambda)$ of the mixed formulation has been established in \cite{haslinger-96} and it holds $\lambda = \sn (u)$. 

\noindent So, the following classical inequality (see \cite{daveiga06}) holds:
\begin{theorem}\label{thm:u-lambda}
Given $s>0$, if the displacement $u$ verifies $u \in H^{3/2+s}(\Omega)^d$, then $\lambda \in H^s(\Gamma_C)$ and it holds:
\begin{eqnarray} \label{ineq:Sp-2 primo}
 \dis \norm{\lambda}_{s,\GC} \leq  \norm{u}_{3/2+s,\Omega}.
\end{eqnarray} 
\end{theorem}

With regards to writing the augmented Lagrange multiplier methods, we use the equivalence between the complementary condition \eqref{eq:contact_cond_mult} and the following equality with $r>0$ a augmented Lagrangian parameter:
\begin{eqnarray} \label{eq:augmultlambda}
\dis\lambda = [\lambda + ru_n]_- , 
\end{eqnarray} 
where $[\cdot]_-$ is the negative part, \textit{i.e.} ${[z]_- = \min ({0,z})}$. This method involved penalizing the multiplier to ensure the contact conditions are verified by the multiplier. \\

\noindent Using this equality, we can express the augmented Lagrangian method (see \cite{al-cu1988,bica-kozia-08,renard-2012b}) as follows
\begin{eqnarray} \label{eq:augmult}
\left\{ 
\begin{array}{l}
\dis a(u,v) +  \int_{\Gamma_C} [\lambda + r u_n]_- v_n \ \textrm{d}\Gamma = L(v)  \quad \forall v, \\
\dis  - \frac{1}{r}   \int_{\Gamma_C} ( \lambda - [\lambda + r u_n]_-)  \mu\ \textrm{d}\Gamma = 0 \quad \forall \mu .
\end{array}
\right.
\end{eqnarray}
Optionally, the second line of the system (\ref{eq:augmult}) can be exploited to replace $\dis [\lambda + ru_n]_-$ in the first line with $\lambda$, we obtain
\begin{eqnarray} \label{eq:v2augmult}
\left\{ 
\begin{array}{l}
\dis a(u,v) +  \int_{\Gamma_C} \lambda v_n \ \textrm{d}\Gamma = L(v)  \quad \forall v, \\
\dis  - \frac{1}{r}   \int_{\Gamma_C} ( \lambda - [\lambda + r u_n]_-)  \mu\ \textrm{d}\Gamma = 0 \quad \forall \mu .
\end{array}
\right.
\end{eqnarray}
We noticed that, the augmented Lagrange multiplier seeks stationary points of the functional: 
\begin{eqnarray} \label{eq:functional}
\begin{array}{l}
\dis  \caL (\varphi, \lambda) := \frac{1}{2} a(u,u) + \frac{1}{2r}\norm{[\lambda + r u_n]_-}^2_{0,\GC} - \frac{1}{2r} \norm{\lambda}^2_{0,\GC}.
\end{array}
\end{eqnarray}
The aim of this paper is to discretize the problem \eqref{eq:augmult} within the isogeometric paradigm, \textit{i.e.} with splines and NURBS. To choose properly the space of Lagrange multipliers properly, we inspire by \cite{brivadis15,fabre-a-priori-iga-17}. In what follows, we introduce NURBS spaces and assumptions together with relevant choices of space pairings. In particular, following \cite{brivadis15,fabre-a-priori-iga-17}, we focus on the definitions of B-Spline displacements of degree $p$ and multiplier spaces of degree $p-2$.

\subsection{NURBS discretisation}
\label{subsec:nurbs_discretisation}

In this section, we give a brief overview on isogeometric analysis providing the notation and concept needed in the next sections. Firstly, we define B-Splines and NURBS in one-dimension. Secondly, we extend these definitions to the multi-dimensional case. Finally, we define the primal and the dual spaces for the contact boundary.

We denote by $Z = \{ \zeta_1, \ldots, \zeta_E\}$ as vector of breakpoints, \textit{i.e.} knots taken without repetition, and $m_j$, the multiplicity of the breakpoint $ \xi_j, \ j=1, \ldots , E$. We define by $p$ the degree of univariate B-Splines and by $\Xi$ an open univariate knot vector, where the first and last entries are repeated $(p + 1)$-times. $\Xi$ is the open knot vector associated to $Z$ where each breakpoint is repeated $m_j$-times, \textit{i.e.} $$\Xi := \{ 0= \xi_1= \cdots =  \xi_{p+1} <  \xi_{p+2} \leq \ldots \leq  \xi_{\eta}< \xi_{\eta+1}= \cdots =  \xi_{\eta+p+1} \}. $$

In what follows, we suppose that $m_1 = m_E = p+1$, while $m_j \leq p-1$, $\forall j = 2,\ldots, E-1$. We define by $\hat{B}^p_i(\zeta)$, $i=1,\ldots,\eta$ the $i$-th univariable B-Spline based on the univariate knot vector $\Xi$ and the degree $p$. We denote by $S^p (\Xi)=  Span\{ \hat{B}^p_i(\zeta), \ i=1,\ldots,\eta \}$. Moreover, for further use we denote by $\tilde{\Xi}$ the sub-vector of $\Xi$ obtained by removing the first and the last knots.

Multivariate B-Splines in dimension $d$ are obtained by tensor product of univariate B-Splines. For any direction $\delta \in \{ 1, \ldots, d\}$, we define by $\eta_\delta$ the number of B-Splines, $\Xi_\delta$ the open knot vector and $Z_\delta$ the breakpoint vector. Then, we define the multivariate knot vector by $\bfX = ( \Xi_1 \times \ldots \times \Xi_d) $ and the multivariate breakpoint vector by $\boldsymbol{Z} = ( Z_1 \times \ldots \times Z_d ) $. We introduce a set of multi-indices $\bfI = \{ \bfi  =(i_1, \ldots, i_d) \mid 1 \leq i_\delta \leq \eta_\delta \}$. We build the multivariate B-Spline functions for each multi-index $\bfi$ by tensorization from the univariate B-Splines, let $\bfz \in \boldsymbol{Z}$ be a parametric coordinate of the generic point: $$\hat{B}^p_{\bfi}(\bfz) = \hat{B}^p_{i_1}(\zeta_1) \ldots \hat{B}^p_{i_d}(\zeta_d) .$$
Let us define the multivariate spline space in the reference domain by (for more details, see \cite{brivadis15}): $$S^p(\bfX) =  Span\{ \hat{B}^p_{\bfi}(\bfz), \ \bfi\in \bfI \}.$$
We define $N^p(\bfX)$ as the NURBS space, spanned by the function $\hat{N}^p_{\bfi}(\bfz)$ with $$\hat{N}^p_{\bfi}(\bfz)  = \frac{\omega_{\bfi}  \hat{B}^p_{\bfi}(\bfz)}{\hat{W}(\bfz)},  $$ where $\{\omega_{\bfi}\}_{\bfi\in \bfI}$ is a set of positive weights and $\dis \hat{W}(\bfz) = \sum_{\bfi\in \bfI} \omega_{\bfi} \hat{B}^p_{\bfi}(\bfz) $ is the weight function and we set $$ N^p(\bfX) =  Span\{ \hat{N}^p_{\bfi}(\bfz) , \ \bfi\in \bfI  \} .$$

In what follows, we will assume that $\Omega$ is obtained as image of $\dis\hat{\Omega} = ]0,1[^d$ through a NURBS mapping $\varphi_0$, \textit{i.e.} $\dis\Omega = \varphi_0(\hat{\Omega})$. Moreover, in order to simplify our presentation, we assume that $\Gamma_C$ is the image of a full face $\dis\hat{f}$ of $\dis\bar{\hat{\Omega}}$, \textit{i.e.} $\dis{\Gamma_C} = \varphi_0(\hat{f})$. We denote by $\dis \varphi_{0,\Gamma_C}$ the restriction of $\dis \varphi_{0}$ to $\dis\hat{f}$. \\

\noindent  A NURBS surface, in d=2, or solid, in d=3, is parameterised by $$ \caC( \bfz) = \sum_{\bfi \in \bfI} C_{\bfi} \hat{N}^p_{\bfi}(\bfz) ,$$
where ${C_{\bfi}}_{ \in \bfI} \in \R^d$, is a set of control point coordinates. The control points are somewhat analogous to nodal points in finite element analysis. The NURBS geometry is defined as the image of the reference domain $\hat{\Omega}$ by $\varphi$, called geometric mapping, $\Omega_t = \varphi(\hat{\Omega})$. \\

We remark that the physical domain $\Omega$ is split into elements by the image of $\boldsymbol{Z}$ through the map $\varphi_0$. We denote such a physical mesh $\caQ_h$ and physical elements in this mesh by $Q$. $\Gamma_C$ inherits a mesh that we denote by $ \restr{\caQ_h}{\GC} $. Elements on this mesh will be denoted as $Q_C$. \\

Finally, we introduce some notations and assumptions on the mesh. 

\noindent \textbf{Assumption 1.} The mapping $\varphi_0$ is considered to be a bi-Lipschitz homeomorphism. Furthermore, for any parametric element ${\hat{Q}}$, $ \dis \restr{\varphi_0}{\bar{\hat{Q}}}$ is in $\caC^\infty(\bar{\hat{Q}})$ and for any physical element ${{Q}}$, $ \dis \restr{\varphi_0^{-1}}{{\bar{Q}}}$ is in $\caC^\infty({\bar{Q}})$.

\noindent Let $h_Q$ be the size of an physical element $Q$, it holds $h_Q = \textrm{diam} (Q)$. In the same way, we define the mesh size for any parametric element. In addition, the Assumption 1 ensures that both size of mesh are equivalent. We denote the maximal mesh size by $\dis h = \max_{Q \in \caQ_h} h_Q$.

\noindent \textbf{Assumption 2.} The mesh $\caQ_h$ is quasi-uniform, \textit{i.e} there exists a constant $\theta$ such that $\dis \frac{h_Q}{h_{Q'}} \leq \theta$ with $Q$ and $Q' \in \caQ_h$.

\section{Discrete spaces and their properties}
\label{sec:discrete_pb}
We will now focus on the definition of spaces on the domain $\Omega$, following the ideas of \cite{brivadis15}. 

\noindent For displacements, we denote by $V^h \subset V$ the space of mapped NURBS of degree $p$ with appropriate homogeneous Dirichlet boundary condition: $$\dis V^h := \{ v^h = \hat{v}^h \circ \varphi^{-1}_0,  \quad \hat{v}^h \in N^p(\bfX)^d \} \cap V .$$
%

\noindent For multipliers, we define the space of B-Splines of degree $p-2$ on the potential contact zone $\dis \Gamma_C = \varphi_{0,\Gamma_C}(\hat{f})$. We denote by $\bfX_{\hat{f}}$ the knot vector defined on $\hat{f}$ and by $\tilde{\bfX}_{\hat{f}}$ the knot vector obtained by removing the first and last value in each knot vector. We define:
$$\dis \Lambda^h := \{ \lambda^h = \hat{\lambda}^h \circ \varphi_{0,\Gamma_C}^{-1}, \quad \hat{\lambda}^h \in S^{p-2}(\tilde{\bfX}_{\hat{f}}) \} .$$
The scalar space $\Lambda^h$ is spanned by mapped B-Splines of the type $\dis \hat{B}^{p-2}_{\bfi}(\bfz)\circ \varphi_{0,\Gamma_C}^{-1}$ for $\bfi$ belonging to a suitable set of indices. In order to reduce our notation, we call $K$ the unrolling of the multi-index $\bfi$,  $K=0 \ldots \caK$ and 
remove super-indices: for $K$ corresponding a given $\bfi$, we set  $ \hat{B}_K (\bfz) = \hat{B}^{p-2}_{\bfi}(\bfz)$, $ {B}_K = \hat{B}_K \circ \varphi_{0,\Gamma_C}^{-1}$  and:
\begin{eqnarray}\label{def:space_mult}
\dis \Lambda^h := Span \{ B_K(x), \quad K=0 \ldots \caK \} . 
\end{eqnarray}
For further use, for $v \in \dis L^2(\GC)$ and for each $K=0 \ldots \caK$, we denote by $\dis (\Pi_\lambda^h \cdot )_K$ the following weighted average of $v$: \begin{eqnarray}\label{def:l2proj_K}
\dis  (\Pi_\lambda^h v)_K = \frac{\dis \intGC v B_K \dG}{\dis \intGC B_K \dG},
\end{eqnarray}
and by $\Pi_\lambda^h$ the global operator such as: \begin{eqnarray}\label{def:l2proj}
\dis \Pi_\lambda^h v = \sum_{K = 0}^\caK (\Pi_\lambda^h  v)_K B_K.
\end{eqnarray}
%
%

We can now define the discrete formulation, as follows:
\begin{eqnarray} \label{eq:augmultdiscrete}
\left\{ 
\begin{array}{l}
\dis a(u^h,v^h) +  \int_{\Gamma_C} [\lambda^h + r (\Pi_\lambda^h u_n^h)]_- \!\   (\Pi_\lambda^h v_n^h) \ \textrm{d}\Gamma = L(v^h)  \quad \forall v^h, \\
\dis  - \frac{1}{r}   \int_{\Gamma_C} ( \lambda^h - [\lambda^h + r (\Pi_\lambda^h u_n^h)]_-)  \mu^h \ \textrm{d}\Gamma = 0 \quad \forall \mu^h .
\end{array}
\right.
\end{eqnarray}
In the following, we define $\dis r =  \frac{r_0}{h}$.

\noindent For any $Q_C \in \restr{\caQ_h}{\GC}$, $\tilde{Q}_C$ denotes the support extension of $Q_C$ (see \cite{daveiga06,daveiga2014}) defined as the image of supports of B-Splines that are not zero on $\hat{Q}_C= \varphi_{0,\Gamma_C}^{-1} (Q_C)$.

\noindent We notice that the operator verifies the following estimate error:
\begin{lemma}\label{lem:op_h1}
Let $\psi \in H^s(\Gamma_C)$ with $0\leq s\leq 1$, the estimate for the local interpolation error reads: 
 \begin{eqnarray} \label{ineq:local disp}
 \dis  \norm{\psi  -\Pi_\lambda^h (\psi )}_{0,\GCQ} \lesssim h^s \norm{\psi}_{s,\GCQtilde} , \qquad \forall \GCQ \in \restr{\caQ_h}{\GC}.
\end{eqnarray}
\end{lemma}
\noindent \pf see the article \cite{fabre-a-priori-iga-17}.

\begin{lemma}\label{lem:ineg_neg}
Let $a,b \in \R$, then we have:
\begin{eqnarray} \label{eq:functional}
\begin{array}{l}
([a]_- - [b]_-)^2 \leq ([a]_- - [b]_-)(a-b),\\
\abs{[a]_- - [b]_-} \leq \abs{a-b}.
\end{array}\nonumber
\end{eqnarray}
\end{lemma}
\noindent \pf Using that for $c \in \R$, $[c]_-^2 = c[c]_- $, it holds:
 \begin{eqnarray} \label{eq:functional}
\begin{array}{lcl}
([a]_- - [b]_-)^2 &=&\dis [a]_-^2 -[a]_-[b]_- -[a]_-[b]_- +[b]_-^2 \\
&\leq&\dis a [a]_- -a [b]_- - b [a]_- +b [b]_- =  ([a]_- - [b]_-)(a-b) .
\end{array} \nonumber
\end{eqnarray}
The second inequality is trivial if $a$ and $b$ have the same sign. If $a>0$ and $b<0$, it holds: $$ \abs{[a]_- - [b]_-} = \abs{b} \leq \abs{a-b}.$$
On the contrary, if $a<0$ and $b>0$, we get: $$ \abs{[a]_- - [b]_-} = \abs{a} \leq \abs{a-b}.$$

In the next section, we prove the existence and the uniqueness of the discrete solution. In the article of Burman's article \cite{burman-16}, they use a Brouwer's fixed point. In order to prove our result \cite{fabre-a-priori-iga-17} in the augmented context, we use the proof using in Nitsche's method \cite{chouly-hild-2013,chouly-hild-renard-2014,franz-overview}, thanks to the hemi-continuity and monotonicity of the operator. First, prove the coercivity property. Then, the existence and uniqueness of the result is deduced from the hemi-continuity of the non-linear operator which corresponds to discrete problem of \eqref{eq:augmultdiscrete}.
Now, we define the following operator $B^h$  from $V^h \times \Lambda^h$ to $V^h \times \Lambda^h$, for all $u^h, v^h \in V^h$ and $\lambda^h, \mu^h \in \Lambda^h$:
$$\dis (B^h(u^h,\lambda^h);(v^h,\mu^h)) :=  a(u^h,v^h) +  \int_{\Gamma_C} [\lambda^h + r (\Pi_\lambda^h u_n^h)]_- \!\  (\Pi_\lambda^h v_n^h)  \ \textrm{d}\Gamma  - \frac{1}{r}   \int_{\Gamma_C} ( \lambda^h - [\lambda^h + r (\Pi_\lambda^h u_n^h)]_-)  \mu^h\ \textrm{d}\Gamma. $$
$$\dis (B^h(u^h,\lambda^h);(v^h,\mu^h)) :=  a(u^h,v^h) +  \frac{1}{r}   \int_{\Gamma_C} [\lambda^h + r (\Pi_\lambda^h u_n^h)]_- (\mu^h +r  (\Pi_\lambda^h v_n^h) ) \ \textrm{d}\Gamma  - \frac{1}{r}   \int_{\Gamma_C} \lambda^h  \mu^h\ \textrm{d}\Gamma. $$
Let us define the following discrete linear operators:
$$ P^h_r : 
 \begin{array}{l}
\dis V^h \times \Lambda^h \rightarrow L^2(\Gamma_C) \\
\dis  (v^h,\mu^h) \mapsto  \mu^h + r (\Pi_\lambda^h v_n^h).
\end{array}\nonumber
, \qquad P_r : 
 \begin{array}{l}
\dis V^h \times \Lambda^h \rightarrow L^2(\Gamma_C) \\
\dis  (v^h,\mu^h) \mapsto  \mu^h + r v_n^h.
\end{array}\nonumber
$$
It holds:
$$\dis (B^h(u^h,\lambda^h);(v^h,\mu^h)) :=  a(u^h,v^h) +  \frac{1}{r}   \int_{\Gamma_C} [ P^h_r(u^h,\lambda^h)]_- P^h_r(v^h,\mu^h) \ \textrm{d}\Gamma  - \frac{1}{r}   \int_{\Gamma_C} \lambda^h  \mu^h\ \textrm{d}\Gamma. $$

First, we need to prove that $B^h$ is coercive.
$$ (B^h(u^h,\lambda^h) - B^h(v^h,\mu^h) ; (u^h,\lambda^h) - (v^h,\mu^h) ) = I + II + III ,$$
with 
\begin{eqnarray}
\begin{array}{l}
\dis I = a(u^h-v^h,u^h-v^h), \\[0.2cm]
\dis II = \frac{1}{r}   \int_{\Gamma_C} \left( [ P^h_r(u^h,\lambda^h)]_- - [P^h_r(v^h,\mu^h)]_-  \right) \left( P^h_r(u^h,\lambda^h) - P^h_r(v^h,\mu^h)  \right)\ \textrm{d}\Gamma , \\[0.4cm]
\dis III =  - \frac{1}{r}   \int_{\Gamma_C} \left( \lambda^h - \mu^h\right)   \left( \lambda^h - \mu^h \right)\ \textrm{d}\Gamma .
\end{array} \nonumber
\end{eqnarray}
We denote by $\alpha$ the ellipticity constant of $a( \cdot, \cdot)$ on $V$, it holds: $$ I \geq \alpha \norm{u^h-v^h}_{1,\Omega}^2 .$$
Using the trace's theorem and Lemma \ref{lem:ineg_neg}, we get:
\begin{eqnarray}
\begin{array}{lcl}
\dis II & = &\dis \frac{1}{r}   \int_{\Gamma_C} \left( [ P^h_r(u^h,\lambda^h)]_- - [P^h_r(v^h,\mu^h)]_-  \right) \left( P^h_r(u^h,\lambda^h) - P^h_r(v^h,\mu^h)  \right)\ \textrm{d}\Gamma \\[0.4cm]
& \geq &\dis \norm{ r^{-\frac{1}{2}}\left([ P^h_r(u^h,\lambda^h)]_- - [P^h_r(v^h,\mu^h)]_- \right)}_{0,\Gamma_C}^2 \\[0.4cm]
\end{array} \nonumber
\end{eqnarray}
And obviously, we have:
\begin{eqnarray}
\begin{array}{l}
\dis III = -\norm{ r^{-\frac{1}{2}}\left( \lambda^h - \mu^h \right)}_{0,\Gamma_C}^2.
\end{array} \nonumber
\end{eqnarray}
We deduce from the previous estimates that:
\begin{eqnarray}
\begin{array}{ll}
\dis  (B^h(u^h,\lambda^h) - B^h(v^h,\mu^h) ; (u^h,\lambda^h) - (v^h,\mu^h) ) \\[0.2cm]
 \dis \qquad ~~~ \geq \alpha \norm{u^h-v^h}_{1,\Omega}^2 +  \norm{ r^{-\frac{1}{2}}\left([ P^h_r(u^h,\lambda^h)]_- - [P^h_r(v^h,\mu^h)]_- \right)}_{0,\Gamma_C}^2  -\norm{ r^{-\frac{1}{2}}\left( \lambda^h - \mu^h \right)}_{0,\Gamma_C}^2 \\[0.4cm]
\end{array} \nonumber
\end{eqnarray}
Thus, if $r_0$ is sufficiently large, we get the coercivity, as follows:
\begin{eqnarray}
\begin{array}{ll}
 (B^h(u^h,\lambda^h) - B^h(v^h,\mu^h) ; (u^h,\lambda^h) - (v^h,\mu^h) ) \geq C \norm{u^h-v^h}_{1,\Omega}^2.
\end{array} \nonumber
\end{eqnarray} \\

Now, we prove the hemi-continuity of $B^h$. Let $s, t \in [0,1]$, $u^h,v^h \in V^h$ and $\lambda^h, \mu^h \in \Lambda^h$, we get:
\begin{eqnarray}
\begin{array}{ll}
 \abs{(B^h(u^h-t v^h,\lambda^h-t\mu^h) - B^h(u^h-s v^h,\lambda^h-s\mu^h) ; (v^h,\mu^h) ) }\\
 \dis  ~~~ \leq \abs{s-t} a(v^h,v^h) + \abs{- \frac{1}{r}   \int_{\Gamma_C} (\lambda^h-t\mu^h)  \mu^h\ \textrm{d}\Gamma  + \frac{1}{r}   \int_{\Gamma_C} (\lambda^h-s\mu^h)  \mu^h\ \textrm{d}\Gamma} \\
\dis  \qquad  + \frac{1}{r}    \int_{\Gamma_C} \abs{ [\lambda^h-t\mu^h + r ((\Pi_\lambda^h u_n^h)-t (\Pi_\lambda^h v_n^h))]_- -  [\lambda^h-s\mu^h + r ((\Pi_\lambda^h u_n^h)-s (\Pi_\lambda^h v_n^h))]_- } \abs{\mu^h +r (\Pi_\lambda^h v_n^h)} \ \textrm{d}\Gamma.
\end{array} \nonumber
\end{eqnarray} 
Using the trace's theorem and Lemma \ref{lem:ineg_neg}, it holds
\begin{eqnarray}
\begin{array}{ll}
 \abs{(B^h(u^h-t v^h,\lambda^h-t\mu^h) - B^h(u^h-s v^h,\lambda^h-s\mu^h) ; (v^h,\mu^h) ) }\\[0.1cm]
 \dis  ~~~ \leq \abs{s-t} a(v^h,v^h) + \abs{s-t}\norm{ r^{-\frac{1}{2}} \mu^h}^2_{0,\Gamma_C}  +  \abs{s-t}  \norm{ r^{-\frac{1}{2}} \left( \mu^h + r (\Pi_\lambda^h v_n^h)\right)}^2_{0,\Gamma_C}.\\[0.3cm]
  \end{array} \nonumber
\end{eqnarray} 
Hence $B^h$ is hemi-continuous.

Let us recall that the following inequalities (see \cite{daveiga06}) are true for the primal and the dual space, before concentrating on the analysis of \eqref{eq:augmultdiscrete}.
\begin{theorem}
Let a given quasi-uniform mesh and let $r, s$ be such that $0 \leq r \leq s \leq p +1$. Then, there exists a constant dependence only on $p, \theta, \varphi_0$ and $\hat{W}$ such that for any $v \in H^s(\Omega)$ there exists an approximation $v^h \in V^h$ such that
\begin{eqnarray} \label{ineq:Np}
\dis \norm{v-v^h}_{r,\Omega} \lesssim h^{s-r} \norm{v}_{s,\Omega}.
\end{eqnarray} 
\end{theorem} 
\noindent We will also make use of the local approximation estimates for splines quasi-interpolants that can be found \textit{e.g.} in \cite{daveiga06,daveiga2014}. 
\begin{lemma}\label{lem:op_mult}
Let $\lambda \in H^s(\GC)$ with $0\leq s \leq p-1$, then there exists a constant depending only on $p,\varphi_0$ and $\theta$, there exists an approximation $\lambda^h \in \Lambda^h$ such that:
\begin{eqnarray} \label{ineq:local}
\dis h^{-1/2} \norm{\lambda-\lambda^h}_{-1/2,\GCQ} +  \norm{\lambda-\lambda^h}_{0,\GCQ} \lesssim h^{s} \norm{\lambda}_{s,\GCQtilde} , \qquad \forall \GCQ \in \restr{\caQ_h}{\GC}.
\end{eqnarray} 
\end{lemma}

  
\section{\textit{A priori} error analysis}
\label{sec:sec2}

In this section, we present an optimal \textit{a priori} error estimate for the unilateral contact problem using augmented Lagrangian method. Our estimates follows the ones for finite elements, provided in Nitsche's context \cite{chouly-hild-2013,chouly-hild-renard-2014,franz-overview}.

For any $p$, we prove our method to be optimal for solutions with regularity up to $5/2$. Thus, optimality for the displacement is obtained for any $p \geq 2$. The cheapest and more convenient method proved optimal corresponds with the choice $p=2$. We also present the choice $p=3$ which produces continuous pressures. Larger values of $p$ may be of interest but, on the other hand, the error bounds remain limited by the regularity of the solution, \textit{i.e.} , up to $C h^{3/2}$. 
In order to prove Theorem \ref{tmp:apriori} which follows, we need a few preparatory Lemmas.

\begin{theorem}\label{thm:pre_apriori} Let $(u,\lambda)$ and $(u^h,\lambda^h)$ be respectively the solution of the augmented Lagrange problem \eqref{eq:mixed_form} and the discrete augmented Lagrange problem \eqref{eq:augmultdiscrete}. Assume that $u \in H^{3/2+\nu}(\Omega)^d$ with $0<\nu<1$. Then, the following error estimate is satisfied:
\begin{eqnarray} \label{eq:pre_apriori}
\begin{array}{l}
\dis \norm{u-u^h}_{1,\Omega}^2 + \norm{ r^{-\frac{1}{2}}\left(\lambda - [ P^h_r(u^h,\lambda^h)]_-  \right) }_{0,\GC}^2 \\[0.2cm]
\dis \qquad \leq C \inf_{v^h \times \mu^h \in V^h \times \Lambda^h} \left( \norm{u-v^h}_{1,\Omega}^2 + \norm{r^{\frac{1}{2}}(u_n-(\Pi_\lambda^h v^h_n))}_{0,\GC}^2 + \norm{r^{-\frac{1}{2}}(\lambda-\mu^h)}_{0,\GC}^2\right).
\end{array}
\end{eqnarray} 
\end{theorem}
\noindent \pf In what follows, we adapt in the IGA and augmented Lagrangian context the proof proved in \cite{fabre-a-priori-iga-17}  .
Using the coercivity of $a(\cdot,\cdot)$, its continuity and Young inequalities, it holds:
\begin{eqnarray}
\begin{array}{ll}
\dis \alpha \norm{u-u^h}_{1,\Omega}^2 & \dis \leq a(u-u^h,u-u^h) = a(u-u^h,u-v^h) + a(u-u^h,v^h-u^h)\\
& \dis \leq C \norm{u-u^h}_{1,\Omega} \norm{u-v^h}_{1,\Omega} +a(u,v^h-u^h) - a(u^h,v^h-u^h)\\[0.4cm]
& \dis \leq \frac{\alpha}{2} \norm{u-u^h}_{1,\Omega}^2 + \frac{C^2}{2 \alpha} \norm{u-v^h}_{1,\Omega}^2 +a(u,v^h-u^h) - a(u^h,v^h-u^h).
\end{array}\nonumber
\end{eqnarray} 
Hence
\begin{eqnarray}
\begin{array}{ll}
\dis \frac{\alpha}{2}  \norm{u-u^h}_{1,\Omega}^2 & \dis \leq  \frac{C^2}{2 \alpha} \norm{u-v^h}_{1,\Omega}^2 +a(u,v^h-u^h) - a(u^h,v^h-u^h).
\end{array}\nonumber
\end{eqnarray} 
Using the definition of operator $P_r$ and equality \eqref{eq:augmultlambda}, we get:
\begin{eqnarray}
\begin{array}{lcl}
IV &=& \dis a(u,v^h-u^h) - a(u^h,v^h-u^h) \\
&  =& \dis \frac{1}{r}  \int_{\Gamma_C} (\lambda - \lambda^h)  (\mu^h- \lambda^h)\ \textrm{d}\Gamma \\[0.2cm]
&&\dis + \frac{1}{r}   \int_{\Gamma_C} \left( [ P^h_r(u^h,\lambda^h)]_- - [P_r(u,\lambda)]_-  \right)P^h_r(v^h - u^h,\mu^h-\lambda^h) \ \textrm{d}\Gamma \\
&=& \dis \norm{r^{-\frac{1}{2}}(\mu^h-\lambda^h)}_{0,\GC}^2 + \frac{1}{r}  \int_{\Gamma_C} (\lambda - \mu^h)  (\mu^h- \lambda^h)\ \textrm{d}\Gamma \\
&&\dis + \frac{1}{r}   \int_{\Gamma_C} \left( [ P^h_r(u^h,\lambda^h)]_- -\lambda  \right) \left( P_r(u,\lambda) - P^h_r(u^h,\lambda^h)  \right)  \ \textrm{d}\Gamma \\[0.4cm]
&&\dis + \frac{1}{r}   \int_{\Gamma_C} \left( [ P^h_r(u^h,\lambda^h)]_- -\lambda  \right) \left( P^h_r(v^h,\mu^h) -  P_r(u,\lambda)  \right)  \ \textrm{d}\Gamma .
\end{array}\nonumber
\end{eqnarray} 
Using Lemma \ref{lem:ineg_neg} and the trace's theorem, it holds: 
\begin{eqnarray}
\begin{array}{lcl}
IV
&\leq& \dis \norm{r^{-\frac{1}{2}}(\mu^h-\lambda^h)}_{0,\GC}^2 + \norm{r^{-\frac{1}{2}}(\lambda-\mu^h)}_{0,\GC} \norm{r^{-\frac{1}{2}}(\lambda^h-\mu^h)}_{0,\GC} \\
&&\dis - \norm{ r^{-\frac{1}{2}}\left(\lambda - [ P^h_r(u^h,\lambda^h)]_-  \right) }_{0,\GC}^2\\[0.4cm]
&& \dis + \norm{ r^{-\frac{1}{2}}\left( \lambda- [ P^h_r(u^h,\lambda^h)]_-  \right) }_{0,\GC} \norm{ r^{-\frac{1}{2}}\left( P_r(u,\lambda) -   P^h_r(v^h,\mu^h)  \right) }_{0,\GC} .
\end{array}\nonumber
\end{eqnarray} 
Using the Young inequalities and Lemma \ref{lem:op_h1} and summing on all elements, it holds:
\begin{eqnarray}
\begin{array}{lcl}
IV &\leq& \dis (1+ \frac{\beta_1}{2}) \norm{r^{-\frac{1}{2}}(\mu^h-\lambda^h)}_{0,\GC}^2 + \frac{1}{2 \beta_1} \norm{r^{-\frac{1}{2}}(\lambda-\mu^h)}_{0,\GC}^2 \\[0.4cm]
&&\dis + (-1 + \frac{1}{2 \beta_2} )\norm{ r^{-\frac{1}{2}}\left(\lambda - [ P^h_r(u^h,\lambda^h)]_-  \right) }_{0,\GC}^2 +  \frac{\beta_2}{2} \norm{ r^{-\frac{1}{2}}\left(P_r(u,\lambda) - P^h_r(v^h,\mu^h)  \right) }_{0,\GC}^2.
\end{array}\nonumber
\end{eqnarray} 
Using the equality \eqref{eq:augmultlambda}, we get:
\begin{eqnarray}
\begin{array}{lcl}
\dis \norm{ r^{-\frac{1}{2}}\left(P_r(u,\lambda) - P^h_r(v^h,\mu^h)  \right) }_{0,\GC}^2 &\leq& \dis 2\norm{ r^{-\frac{1}{2}}\left(\lambda -\mu^h  \right) }_{0,\GC}^2 + 2\norm{ r^{\frac{1}{2}}\left(u_n- (\Pi_\lambda^h v^h_n)  \right) }_{0,\GC}^2 .
\end{array}\nonumber
\end{eqnarray} 
It holds the following inequality:
\begin{eqnarray}
\begin{array}{lcl}
\dis \norm{r^{-\frac{1}{2}}(\mu^h-\lambda^h)}_{0,\GC}^2 \leq 2\norm{r^{-\frac{1}{2}}(\mu^h-\lambda)}_{0,\GC}^2 +  2\norm{r^{-\frac{1}{2}}(\lambda-\lambda^h)}_{0,\GC}^2 .
\end{array}\nonumber
\end{eqnarray} 
If $\beta_2$ is chosen sufficiently large such that $$\dis -1 + \frac{1}{2 \beta_2}   \leq -\frac{1}{2}  .$$
And if $r_0$ is sufficiently large, this ends the proof of Theorem \ref{thm:pre_apriori}.
\cqfd

\begin{theorem}\label{tmp:apriori} Let $(u,\lambda)$ and $(u^h,\lambda^h)$ be respectively the solution of the mixed problem \eqref{eq:mixed_form} and the discrete augmented Lagrange problem \eqref{eq:augmultdiscrete}. Assume that $u \in H^{3/2+\nu}(\Omega)^d$ with $0<\nu<1$. Then, the following error estimate is satisfied:
\begin{eqnarray} \label{eq:apriori}
\dis \norm{u-u^h}_{1,\Omega}^2 + \norm{ r^{-\frac{1}{2}}\left(\lambda - [ P^h_r(u^h,\lambda^h)]_-  \right) }_{0,\GC}^2 \lesssim h^{1+2\nu} \norm{u}^2_{3/2 + \nu,\Omega}.
\end{eqnarray} 
\end{theorem}
\noindent \pf Now we can establish the inequality \eqref{eq:apriori}, for $r = r_0 / h$, if we replace $(v^h,\mu^h)$ by $(u^h,\lambda^h)$ and using Theorem \ref{ineq:Np}, Lemmas \ref{lem:op_h1} and \ref{lem:op_mult} and summing on all elements, in the inequality \eqref{eq:pre_apriori}. 
\cqfd

%
\section{Numerical Study} 
\label{sec:sec3}
{In this section, we perform {a} numerical validation for the method we propose in small as well as in large deformation {frameworks}, \textit{i.e.}, also beyond the theory developed in previous Sections. Due to the intrinsic lack of regularity of contact solutions, we restrict ourselves to the case $p=2$,  {for which the $N_2 / S_0$ method is tested}, and $p=3$,  {for which the $N_3 / S_1$ method is tested}.}
The suite of benchmarks reproduces the classical Hertz contact problem \cite{Hertz1882,johnsonbook}:
Sections \ref{subsec:Hertz problems 2d} and \ref{subsec:Hertz problems 3d} analyse
the two and three-dimensional cases for a small deformation setting, whereas
Section \ref{subsec:large deformation} considers the large deformation problem in 2D.
The examples were performed using an in-house code based on the igatools library (see \cite{Pauletti2015} for further details).

{In the following example, to prevent that the {contact zone} is empty, we considered, only for the initial gap, that {there exists contact} if the $\dis g_n\leq 10^{-9}$.}\\

\subsection{Two-dimensional Hertz problem}
\label{subsec:Hertz problems 2d}
\begin{figure}[!ht]
  \centering
  \subfigure[Stress magnitude distribution for the undeformed mesh with $N_2 / S_0$ method.]{\label{fig:mesh1_a} \includegraphics[width= 8.5cm]{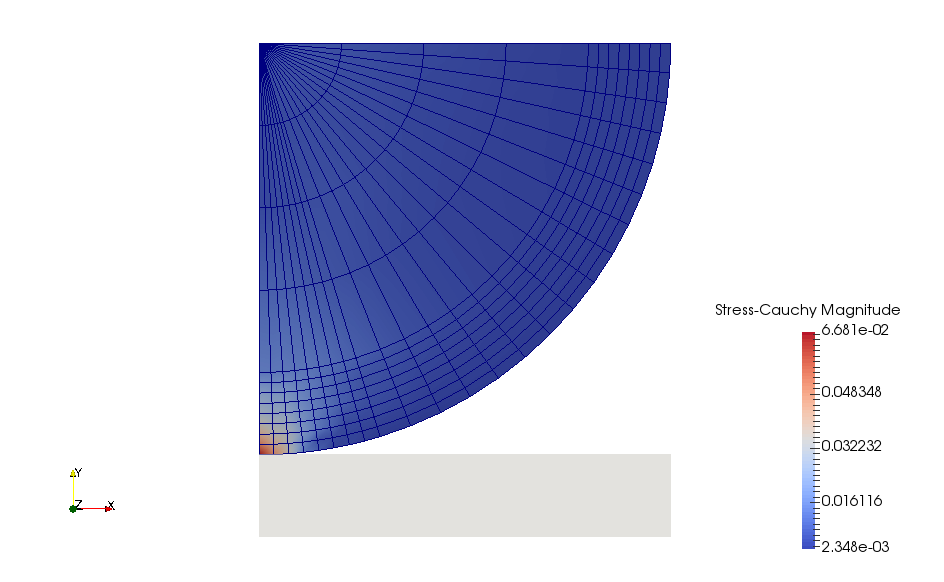}}\hfill
  \subfigure[Analytical and numerical contact pressure with $N_2 / S_0$ and $N_3 / S_1$ methods.]{\label{fig:mesh1_b} \includegraphics[width= 7cm]{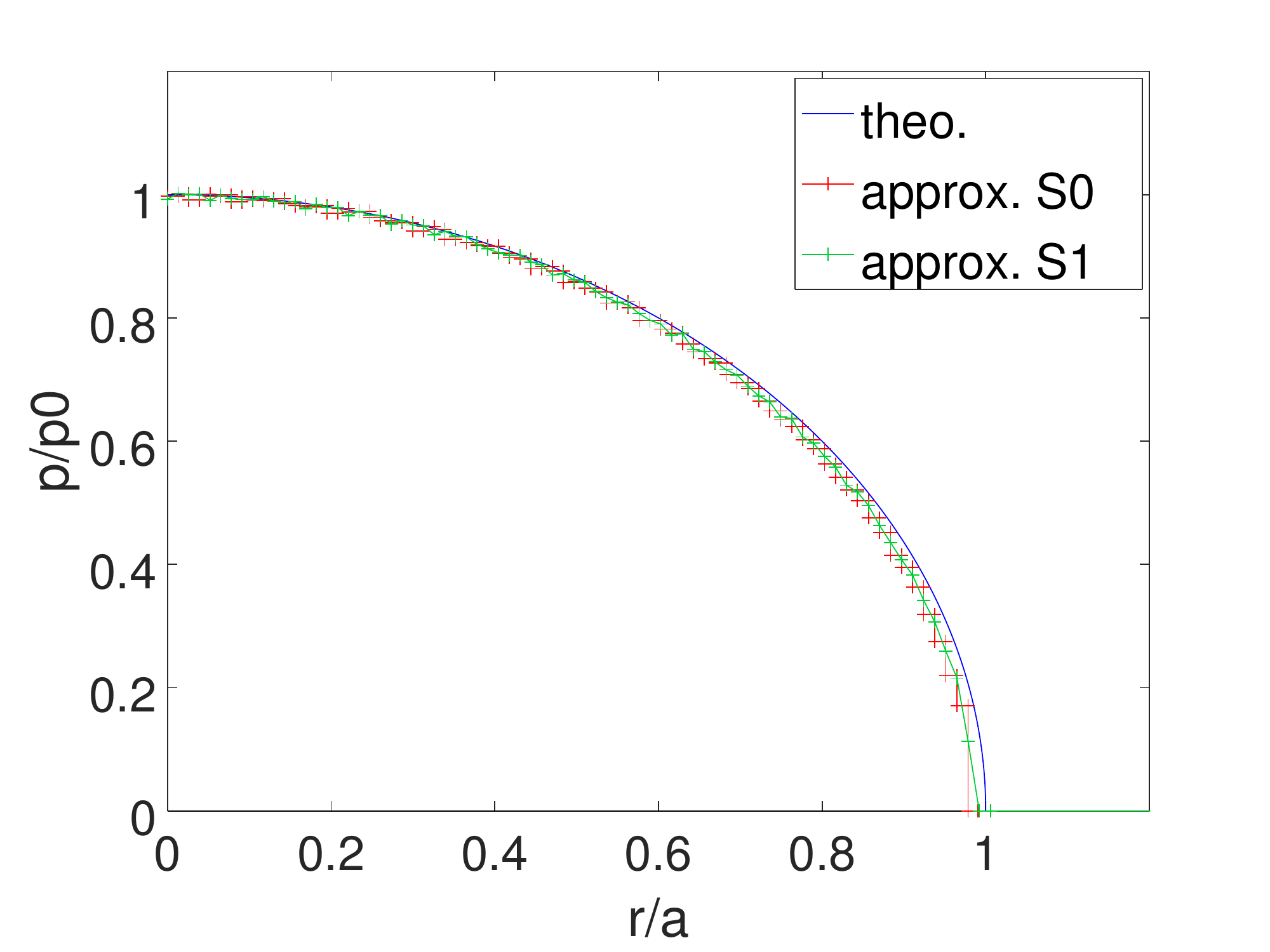}}
  \caption{2D Hertz contact problem  for an applied pressure $P=0.003$.}
  \label{fig:mesh1}
\end{figure}

The first example included in this section analyses the two-dimensional frictionless Hertz contact problem considering small elastic deformations.
It consists in an infinitely long half cylinder body with radius $R=1$, that it is deformable and whose
material is linear elastic, with Young's modulus $E = 1$ and Poisson's ratio $\nu = 0.3$.
A uniform pressure $P=0.003$ is applied on the top face of the cylinder while the curved surface contacts against a horizontal rigid plane
(see Figure \ref{fig:mesh1_a}).
Taking into account the test symmetry and the ideally infinite length of the cylinder,
the problem is modelled as 2D quarter of disc with proper boundary conditions.

Under the hypothesis that the contact area is small compared to the cylinder dimensions, the Hertz's analytical solution (see \cite{Hertz1882,johnsonbook})
predicts that the contact region is an infinitely long band whose width is $2a$, being $a = \sqrt{8R^2P(1- \nu^2)/\pi E}$.
Thus, the normal pressure, {that} follows an elliptical distribution along the width direction $r$, is $p(r)=p_0\sqrt{1-r^2/a^2}$,
where the maximum pressure, at the central line of the band ($r=0$), is $p_0 = 4 RP/ \pi a$.
For the geometrical, material and load data chosen in this numerical test,
the characteristic values of the solution are $a=0.083378$ and $p_0 = 0.045812$.
Notice that, as required by Hertz's theory hypotheses, $a$ is sufficiently small compared to $R$.

It is important to remark that, despite the fact that Hertz's theory provides a
full description of the contact area and the normal contact pressure
in the region, it does not describe analytically the deformation of the whole elastic domain.
Therefore, for all the test cases hereinafter, the $L^2$ {error norm} and $H^1$ {error semi-norm} of
the displacement obtained numerically are computed taking a more refined solution
as a reference. For this bidimensional test case, the mesh size
of the refined solution $h_{ref}$ is such that, for all the discretizations,
$4 h_{ref} \leq h$, where $h$ is the size of the mesh considered.
Additionally, as it is shown in Figure \ref{fig:mesh1_a}, the mesh is finer in the
vicinity of the potential contact zone. The knot vector values
are defined such that $80\%$ of the knot spans are located within $10\%$ of the total length of the knot vector.

In particular, the analysis of this example focuses on the effect of the
interpolation order on the quality of contact stress distribution.
Thus, in Figure \ref{fig:mesh1_b} we  compare the pressure reference solution {with the Lagrange multiplier values computed at the control points for $S0$ elements, \textit{i.e.} its constant values, and for $S1$ elements.}
The dimensionless contact pressure $p/p_0$ is plotted respect to the normalized
coordinate $r/a$.
{For both, the results are very good: the maximum pressure computed
and the pressure distribution, even across the boundary of the contact region (on the contact and non contact zones), are close to the analytical solution.}

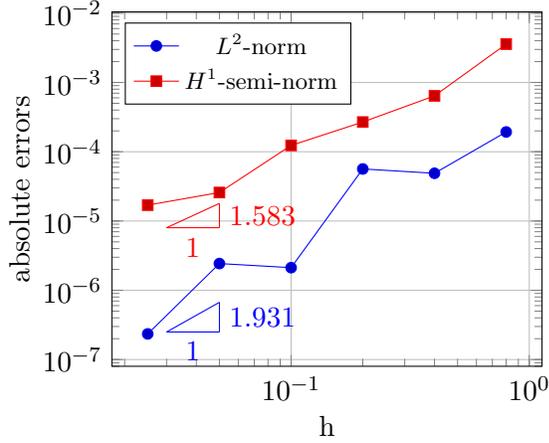
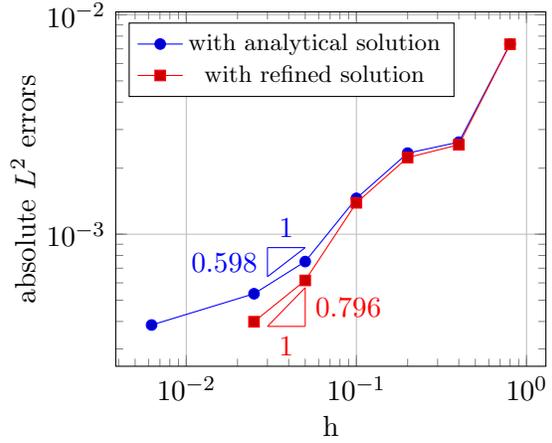
\begin{figure}[!ht]
 \centering
 \subfigure[Displacement error.]{\label{fig:disp 0point003_a}
\begin{tikzpicture}
\begin{loglogaxis}[width=7.3cm,xlabel={h},ymin=0.8e-07,ylabel={absolute errors},legend pos={north west},grid=major]

\addplot table[x=h,y=L2_abs] {figure_latex_figure_2d_hertz_p_0_003_N2_S0_data_disp_N2_S0_p_0point003.res};
\addplot table[x=h,y=H1_abs] {figure_latex_figure_2d_hertz_p_0_003_N2_S0_data_disp_N2_S0_p_0point003.res};

\slopeTriangleAbove{0.03}{0.05}{2.5e-07}{6.7039e-07}{1.931}{blue}; 
\slopeTriangleAbove{0.03}{0.05}{8e-06}{1.796e-05}{1.583}{red}; 
\legend{ \footnotesize $L^2$-norm, \footnotesize $H^1$-{semi-}norm}
\end{loglogaxis}
\end{tikzpicture}
} \hfill
 \subfigure[Lagrange multipliers error.]{\label{fig:disp 0point003_b}

\begin{tikzpicture}
\begin{loglogaxis}[width=7.3cm,xlabel={h},ymin=2.5e-04,ylabel={absolute $L^2$ errors},legend pos={north west},grid=major]

\addplot table[x=h_mult_ana,y=L2_mult_abs_ana] {figure_latex_figure_2d_hertz_p_0_003_N2_S0_data_mult_analytical_N2_S0_p_0point003.res};
\addplot table[x=h_mult_ref,y=L2_mult_abs_ref] {figure_latex_figure_2d_hertz_p_0_003_N2_S0_data_mult_refined_N2_S0_p_0point003.res};

\slopeTriangleBelow{0.03}{0.05}{6.4e-04}{8.6865e-04}{0.598}{blue}; 
\slopeTriangleAbove{0.03}{0.05}{3.8e-04}{5.70657e-04}{0.796}{red}; 
\legend{ \footnotesize with analytical solution,\footnotesize with refined solution}
\end{loglogaxis}
\end{tikzpicture}
}
 \caption{2D Hertz contact problem with $N_2 / S_0$ method for an applied pressure $P=0.003$.
 Absolute displacement errors in $L^2$-norm and $H^1$-{semi-}norm and Lagrange multipliers error in $L^2$-norm, respect to analytical and refined numerical solutions.} 
 \label{fig:disp 0point003}
\end{figure}

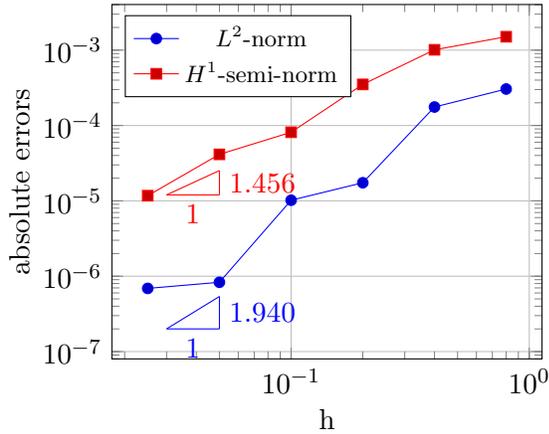
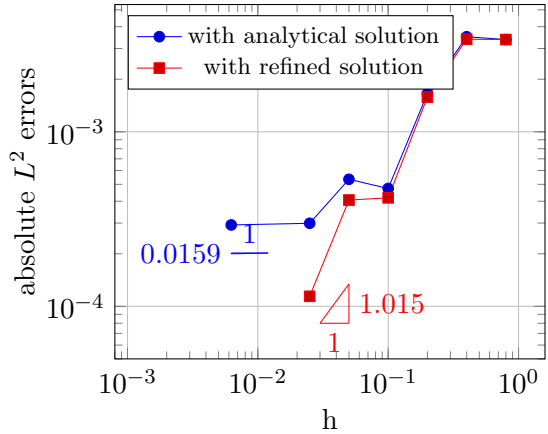
\begin{figure}[!ht]
 \centering
 \subfigure[Displacement error.]{\label{fig:sec3 disp 0point003 2} 
\begin{tikzpicture}
\begin{loglogaxis}[width=7.3cm,xlabel={h},ymin=0.8e-07,ylabel={absolute errors},legend pos={north west},grid=major]

\addplot table[x=h,y=L2_abs] {figure_latex_figure_2d_hertz_p_0_003_N3_S1_data_disp_N3_S1_p_0point003.res};
\addplot table[x=h,y=H1_abs] {figure_latex_figure_2d_hertz_p_0_003_N3_S1_data_disp_N3_S1_p_0point003.res};

\slopeTriangleAbove{0.03}{0.05}{2e-07}{5.3879e-07}{1.940}{blue}; 
\slopeTriangleAbove{0.03}{0.05}{12e-06}{2.5246e-05}{1.456}{red}; 
\legend{ \footnotesize $L^2$-norm, \footnotesize $H^1$-{semi-}norm}
\end{loglogaxis}
\end{tikzpicture}
} \hfill
 \subfigure[Lagrange multipliers error.]{\label{fig:sec3 mult 0point003 2} 

\begin{tikzpicture}
\begin{loglogaxis}[width=7.3cm,xlabel={h},xmin=0.8e-03,ymin=0.5e-04,ylabel={absolute $L^2$ errors},legend pos={north west},grid=major]

\addplot table[x=h_mult_ana,y=L2_mult_abs_ana] {figure_latex_figure_2d_hertz_p_0_003_N3_S1_data_mult_analytical_N3_S1_p_0point003.res};
\addplot table[x=h_mult_ref,y=L2_mult_abs_ref] {figure_latex_figure_2d_hertz_p_0_003_N3_S1_data_mult_refined_N3_S1_p_0point003.res};

\slopeTriangleBelow{0.0062500}{0.012}{2e-04}{2.0209e-04}{0.0159}{blue}; 
\slopeTriangleAbove{0.03}{0.05}{0.8e-04}{1.3436e-04}{1.015}{red}; 
\legend{ \footnotesize with analytical solution,\footnotesize with refined solution}
\end{loglogaxis}
\end{tikzpicture}
}
 \caption{2D Hertz contact problem with $N_3 / S_1$ method for an applied pressure $P=0.003$.
 Absolute displacement errors in $L^2$-norm and $H^1$-{semi-}norm and Lagrange multipliers error in $L^2$-norm, respect to analytical and refined numerical solutions.} 
 \label{fig:disp 0point003 2}
\end{figure}

In Figure \ref{fig:disp 0point003_a} and respectively \ref{fig:sec3 disp 0point003 2}, absolute errors in $L^2${-norm} and $H^1${-semi-norm}
for the $N_2 / S_0$ and respectively $N_3 / S_1$ choice are shown. As expected, optimal convergence are obtained
for the displacement error in the $H^1$-{semi-}norm: the convergence rate is close to the expected $3/2$ value.
Nevertheless, the $L^2$-norm of the displacement error presents suboptimal convergence (close to $2$), but
according to Aubin-Nitsche's lemma in the linear case, the expected convergence rate is {lower than} $5/2$.
On the other hand, in Figure \ref{fig:disp 0point003_b} the 
$L^2$-norm of the Lagrange multipliers error is presented, the expected convergence rate is $1$.
Whereas a convergence rate close to $0.6$ is achieved when {we compare the numerical solution and the Hertz's analytical solution}, and close to $0.8$ is achieved when {we compare the numerical solution and the refined numerical solution}.
In Figure \ref{fig:sec3 mult 0point003 2}, we seem reach a ceiling when we compare the numerical solution and the Hertz's analytical solution and a convergence rate close to $1$ is achieved when we compare the numerical solution and the refined numerical solution.

As a second example, we present the same test case but with significantly higher pressure applied
$P = 0.01$. Under these load conditions, the contact area is wider ($a =  0.15223$) and the
contact pressure higher ($p_0 = 0.083641$).
It can be considered that the ratio $a/R$ no longer satisfies the hypotheses of Hertz's theory.

\begin{figure}[!ht]
  \centering
  \subfigure[Stress magnitude distribution for the undeformed mesh with $N_2 / S_0$ method.]{\label{fig:mesh2_a} \includegraphics[width= 8.5cm]{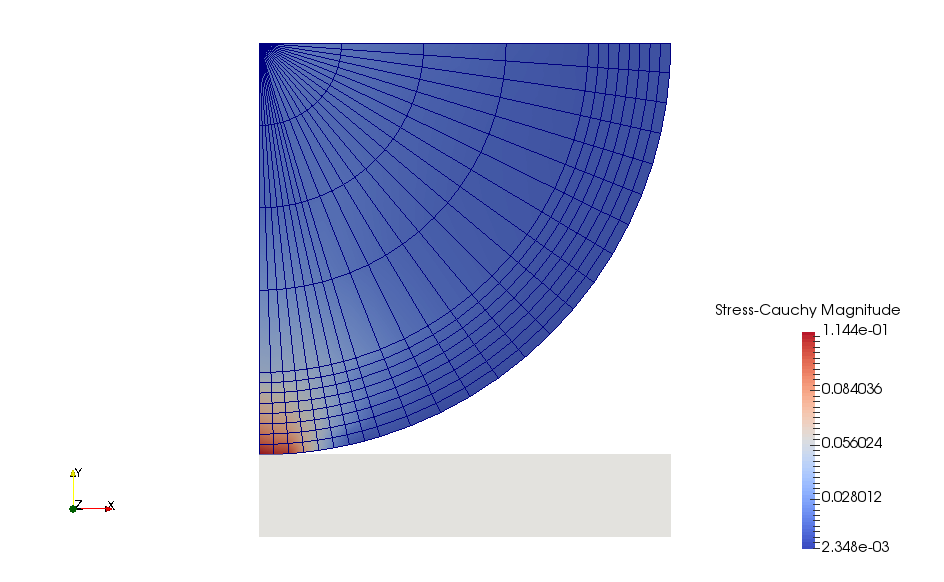}}\hfill
  \subfigure[Analytical and numerical contact pressure with $N_2 / S_0$ and $N_3 / S_1$ methods]{\label{fig:mesh2_b} \includegraphics[width= 7cm]{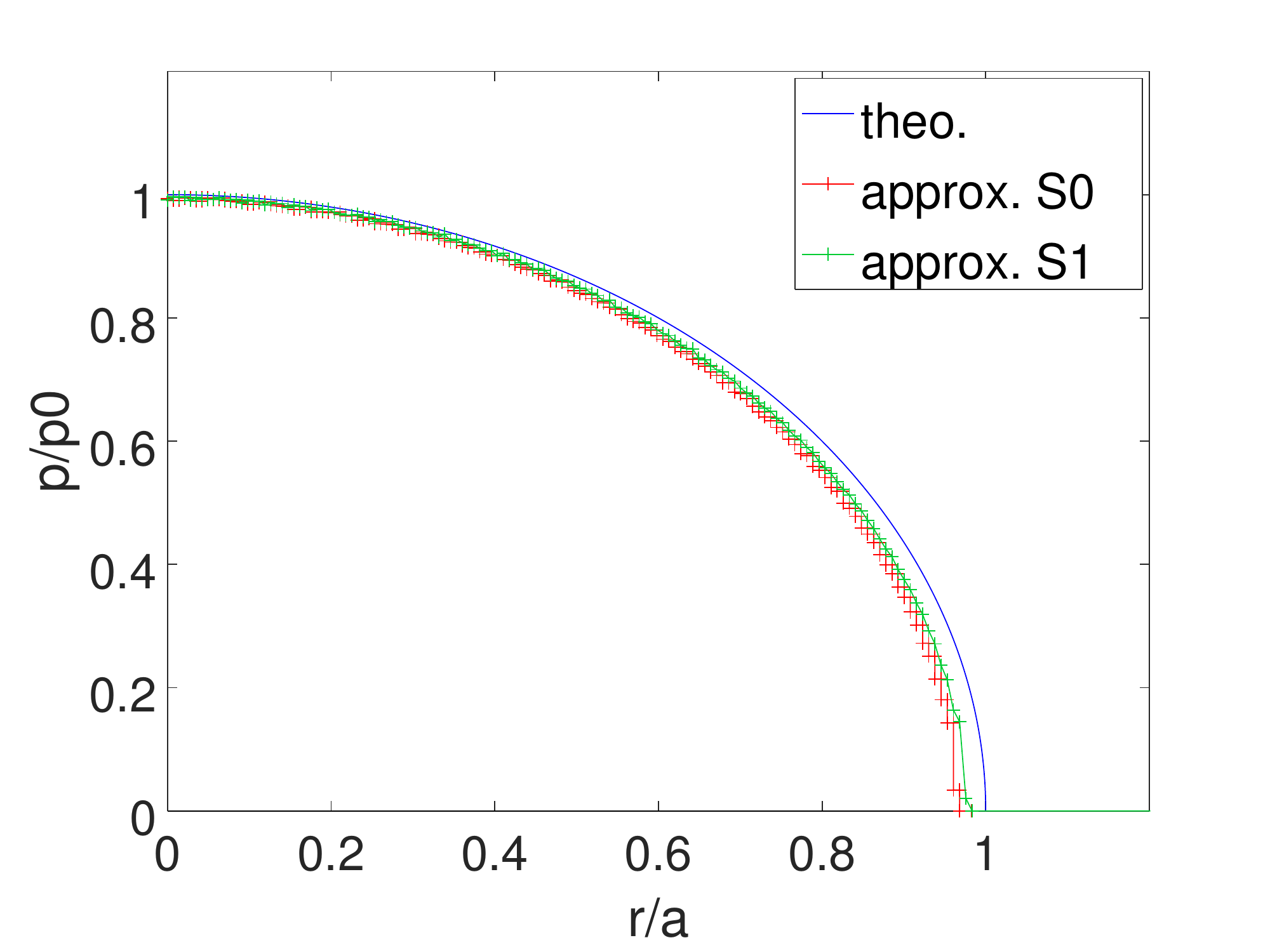}}
  \caption{2D Hertz contact problem  for an applied pressure $P=0.01$.}
  \label{fig:mesh2}
\end{figure}

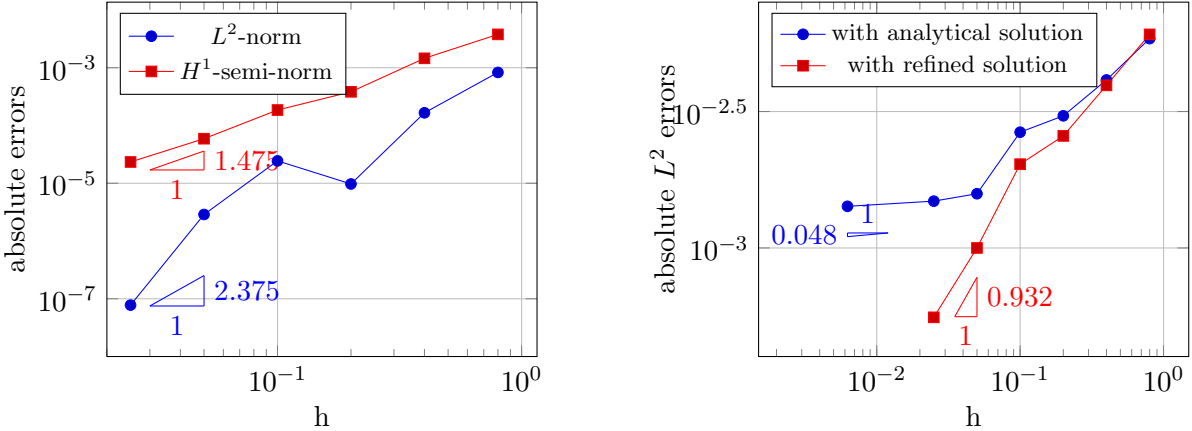
\begin{figure}[!ht]
 \centering
 \subfigure[Displacement error.]{\label{fig:disp 0point01_a} 
 
 \begin{tikzpicture}
\begin{loglogaxis}[width=7.3cm,ymin=1e-08,xmin=2e-02,xlabel={h},ylabel={absolute errors},legend pos={north west},grid=major]

\addplot table[x=h,y=L2_abs] {figure_latex_figure_2d_hertz_p_0_01_N2_S0_data_disp_N2_S0_p_0point01.res};
\addplot table[x=h,y=H1_abs] {figure_latex_figure_2d_hertz_p_0_01_N2_S0_data_disp_N2_S0_p_0point01.res};

\slopeTriangleAbove{0.03}{0.05}{7.5e-08}{2.52320e-07}{2.375}{blue}; 
\slopeTriangleAbove{0.03}{0.05}{1.7e-05}{3.6114e-05}{1.475}{red}; 
\legend{ \footnotesize $L^2$-norm, \footnotesize $H^1$-{semi-}norm}

\end{loglogaxis}
\end{tikzpicture}
 
 }\hfill
 \subfigure[Lagrange multipliers error.]{\label{fig:disp 0point01_b}
 
\begin{tikzpicture}
\begin{loglogaxis}[every axis y label/.style={at={(ticklabel cs:0.5)},rotate=90,anchor=center},width=7.3cm,xlabel={h},ymin=4e-04,ylabel={absolute $L^2$ errors},legend pos={north west},xmin=1.5e-3,grid=major]

\addplot table[x=h_mult_ana,y=L2_mult_abs_ana] {figure_latex_figure_2d_hertz_p_0_01_N2_S0_data_mult_analytical_N2_S0_p_0point01.res};
\addplot table[x=h_mult_ref,y=L2_mult_abs_ref] {figure_latex_figure_2d_hertz_p_0_01_N2_S0_data_mult_refined_N2_S0_p_0point01.res};

\slopeTriangleBelow{0.0062500}{0.012}{0.0011}{0.00113498767}{0.048}{blue}; 
\slopeTriangleAbove{0.035}{0.05}{5.6e-04}{7.8083e-04}{0.932}{red}; 
\legend{ \footnotesize with analytical solution,\footnotesize with refined solution}

\end{loglogaxis}
\end{tikzpicture}  
  
  }
\caption{2D Hertz contact problem with $N_2 / S_0$ method for an applied pressure $P=0.01$.
 Absolute displacement errors in $L^2$-norm and $H^1$-{semi-}norm and Lagrange multipliers error in $L^2$-norm, respect to analytical and refined numerical solutions.} 
 \label{fig:disp 0point01}
\end{figure}

\begin{figure}[!ht]
 \centering
 \subfigure[Displacement error.]{\label{fig:sec3 disp 0point01 2} 
 
 \begin{tikzpicture}
\begin{loglogaxis}[width=7.3cm,ymin=1e-07,xmin=1e-02,xlabel={h},ylabel={absolute errors},legend pos={north west},grid=major]

\addplot table[x=h,y=L2_abs] {figure_latex_figure_2d_hertz_p_0_01_N3_S1_data_disp_N3_S1_p_0point01.res};
\addplot table[x=h,y=H1_abs] {figure_latex_figure_2d_hertz_p_0_01_N3_S1_data_disp_N3_S1_p_0point01.res};

\slopeTriangleBelow{0.03}{0.05}{2e-06}{6.9203e-06}{2.430}{blue}; 
\slopeTriangleBelow{0.03}{0.05}{8e-05}{1.8689e-04}{1.661}{red}; 
\legend{ \footnotesize $L^2$-norm, \footnotesize $H^1$-{semi-}norm}

\end{loglogaxis}
\end{tikzpicture}
 
 }\hfill
 \subfigure[Lagrange multipliers error.]{\label{fig:sec3 mult 0point01 2}
 
\begin{tikzpicture}
\begin{loglogaxis}[every axis y label/.style={at={(ticklabel cs:0.5)},rotate=90,anchor=center},width=7.3cm,xlabel={h},ymin=1e-04,ylabel={absolute $L^2$ errors},legend pos={north west},xmin=1.5e-3,grid=major]

\addplot table[x=h_mult_ana,y=L2_mult_abs_ana] {figure_latex_figure_2d_hertz_p_0_01_N3_S1_data_mult_analytical_N3_S1_p_0point01.res};
\addplot table[x=h_mult_ref,y=L2_mult_abs_ref] {figure_latex_figure_2d_hertz_p_0_01_N3_S1_data_mult_refined_N3_S1_p_0point01.res};

\slopeTriangleBelow{0.0062500}{0.012}{0.002}{  0.0020798}{0.060}{blue}; 
\slopeTriangleBelow{0.03}{0.05}{0.00025}{4.8443e-04}{1.295}{red}; 
\legend{ \footnotesize with analytical solution,\footnotesize with refined solution}

\end{loglogaxis}
\end{tikzpicture}  
  
  }
\caption{2D Hertz contact problem with $N_3 / S_1$ method for an applied pressure $P=0.01$.
 Absolute displacement errors in $L^2$-norm and $H^1$-{semi-}norm and Lagrange multipliers error in $L^2$-norm, respect to analytical and refined numerical solutions.} 
 \label{fig:disp 0point01 2}
\end{figure}
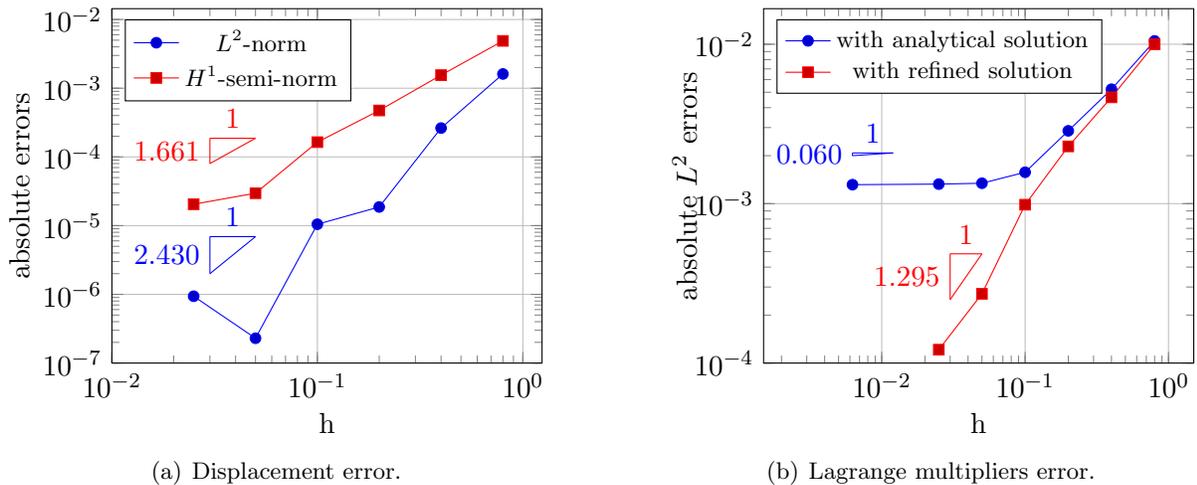

In the same way as before, Figure \ref{fig:mesh2} shows the stress tensor magnitude
and computed contact pressure.
Figure \ref{fig:disp 0point01_a} and respectively Figure \ref{fig:sec3 disp 0point01 2} show the displacement absolute error in $L^2$-norm and $H^1$-{semi-}norm
for $N_2 / S_0$ method and respectively for $N_3 / S_1$ method. As expected, optimal convergence is obtained in the $H^1$-{semi-}norm,
(the convergence rate is close to $1.5$) and, {while, for the $L^2$-norm we obtain a better rate (as expected by the Aubin-Nische's lemma) which can hardly be estimated precisely from the graph.} 
On the other hand, in Figure \ref{fig:disp 0point01_b} and \ref{fig:sec3 mult 0point01 2} it can be seen that the $L^2$-norm of the error of the Lagrange multipliers
evidences a suboptimal behaviour: the error, that initially decreases, remains constant for smaller values of $h$.
It may due to the choice of an excessively large normal pressure:
the approximated solution converges, but not to the analytical solution, that is no longer valid.
Indeed, when compared to a refined numerical solution (Figure \ref{fig:disp 0point01_b} and \ref{fig:sec3 mult 0point01 2}),
the computed Lagrange multipliers solution converges optimally for $N_2 / S_0$-method and sub-optimally for $N_3 / S_1$-method. As it was pointed out above,
for these examples the displacement solution error is computed respect to a more refined numerical solution, therefore, this effect
is not present in displacement results. 

\subsection{Three-dimensional Hertz problem}
\label{subsec:Hertz problems 3d}
\begin{figure}[!ht]
  \centering
  \subfigure[Stress magnitude distribution for the undeformed mesh.]{\label{fig:3d p = 0.0005_a} \includegraphics[width= 8.5cm]{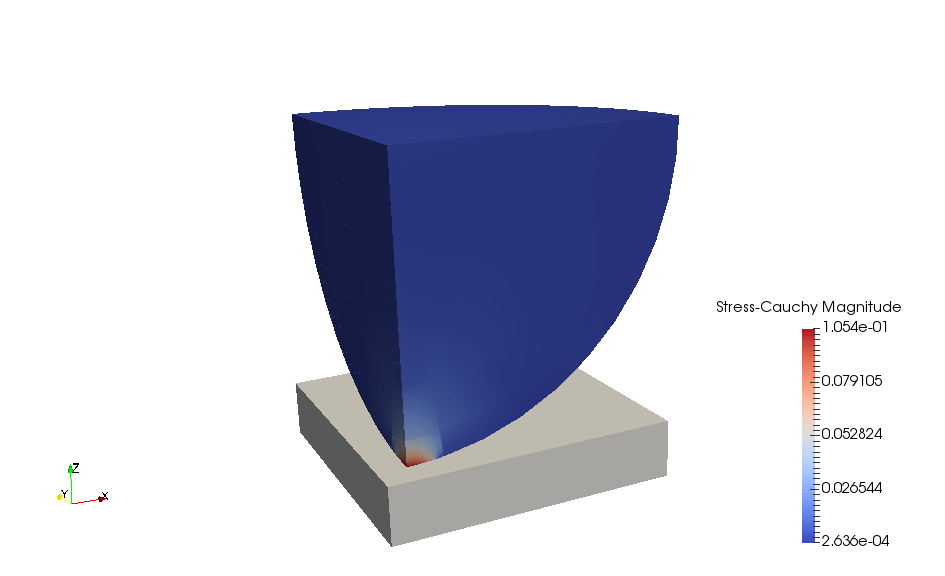}}\hfill
  \subfigure[Analytical and numerical contact pressure for $h = 0.15$.]{\label{fig:3d p = 0.0005_b} \includegraphics[width= 7cm]{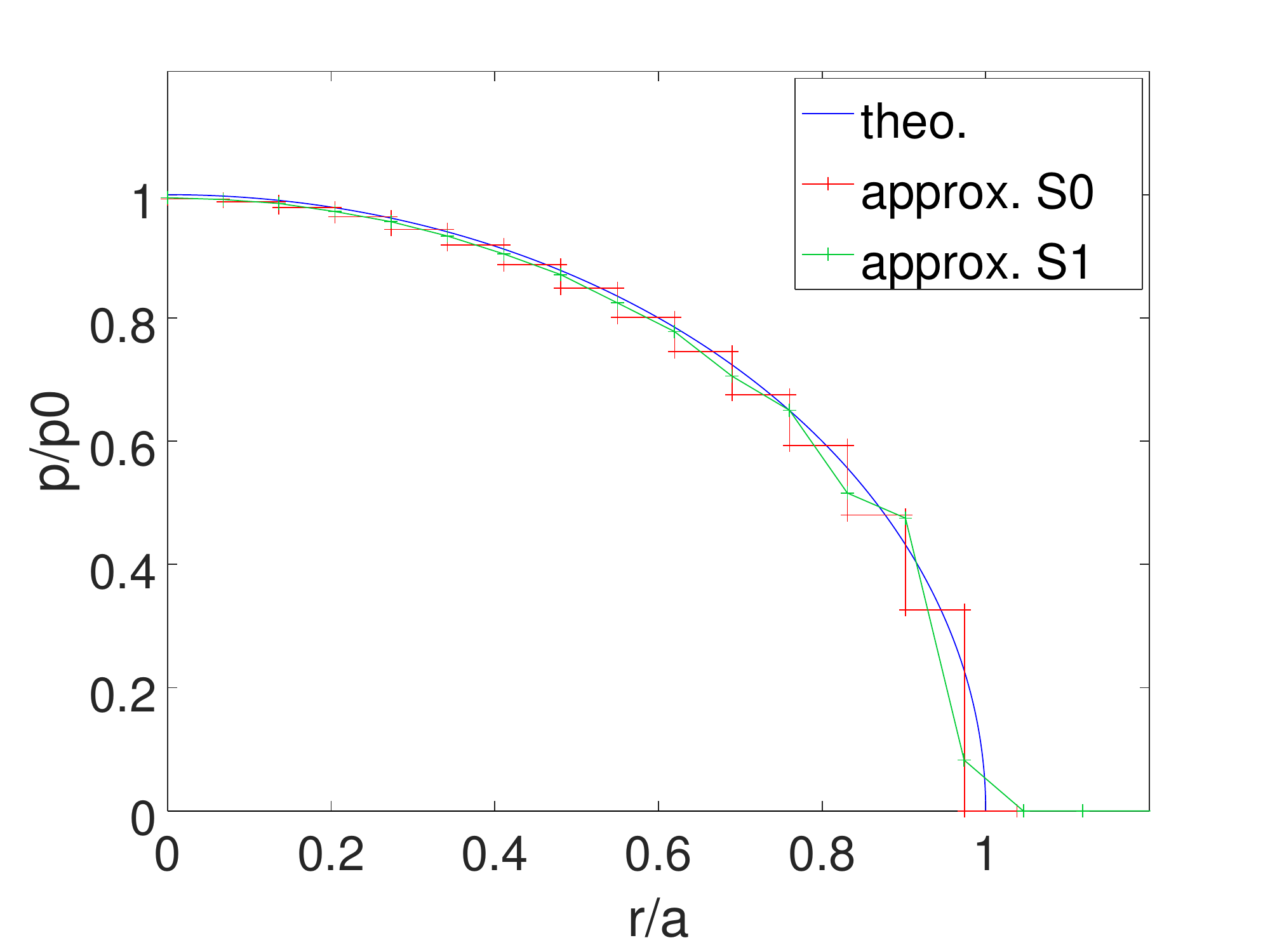}}
  \caption{3D Hertz contact problem with $N_2 / S_0$ method for a higher pressure ($P=5\cdot10^{-4}$).}
  \label{fig:3d p = 0.0005}
\end{figure}

In this section, the three-dimensional frictionless Hertz problem is studied.
It consists in a hemispherical elastic body with radius $R$ that contacts against a horizontal rigid plane as a consequence of
an uniform pressure $P$ applied on the top face (see Figure \ref{fig:3d p = 0.0005_a}).
Hertz's theory predicts that the contact region is a circle of radius $a = (3R^3P(1- \nu^2)/4 E)^{1/3}$
and the contact pressure follows a hemispherical distribution
$p(r)=p_0 \sqrt{1-r^2/a^2}$, with $p_0 = 3R^2P/ 2 a^2$, being $r$ the distance to the centre of the circle
(see\cite{Hertz1882,johnsonbook}).
In this case, for the chosen values $R=1$, $E=1$, $\nu=0.3$ and $P=5\cdot10^{-4}$, the contact radius is $a = 0.10235$ and the maximum pressure $p_0 = 0.0716$.
As in the two-dimensional case, Hertz's theory relies on the hypothesis that $a$ is small compared to $R$ and the deformations are small.

\begin{figure}[!ht]
  \centering
  \subfigure[$h=0.4$.]{\label{fig:3d cp 0.0005_a} \includegraphics[width= 7cm]{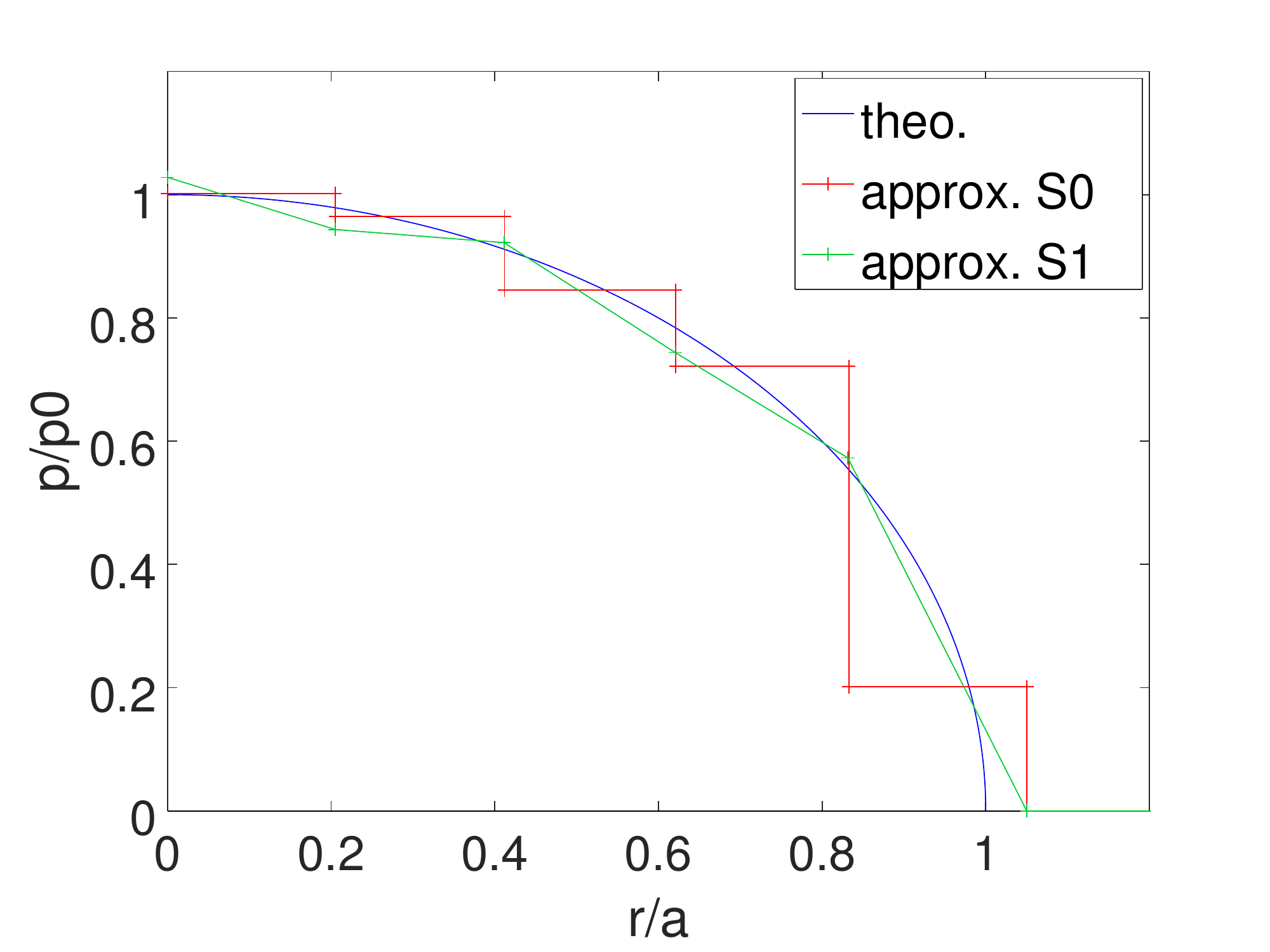}}\hfill
  \subfigure[$h=0.2$.]{\label{fig:3d cp 0.0005_b} \includegraphics[width= 7cm]{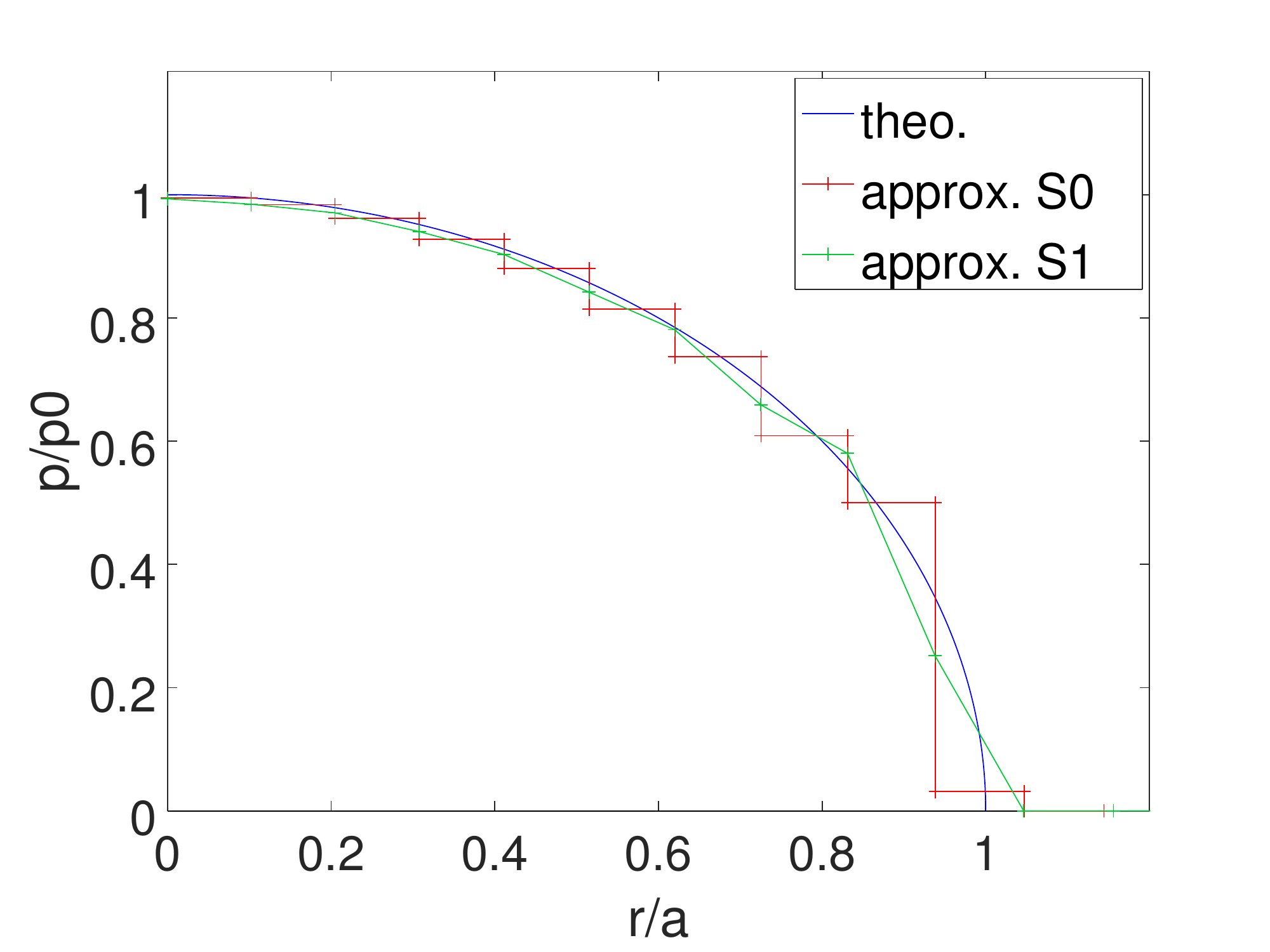}}
  \caption{3D Hertz contact problem with $N_2 / S_0$ method for an applied pressure $P=5\cdot10^{-4}$. Contact pressure solution at control points.}
  \label{fig:3d cp 0.0005}
\end{figure}

Considering the problem axial symmetry, the test is reproduced using an octant of sphere with proper boundary conditions.
Figure \ref{fig:3d cp 0.0005_a} shows the problem setup and the magnitude of the computed stresses.
As in the 2D case, in order to achieve more accurate results in the contact region, the mesh is refined in the vicinity of the potential contact zone.
The knot vectors are defined such as $75\%$ of the elements are located within $10\%$ of the total length of the knot vector.
On the other hand, in Figure \ref{fig:3d cp 0.0005} the contact pressure is shown at control points for mesh sizes $h=0.4$ and $h=0.2$. As it can be
appreciated, good agreement between the analytical and computed pressure is obtained in all cases.
Due to the coarse reference mesh, it is not possible to present good curves of convergence and show the asymptotical behaviour.

\subsection{Two-dimensional Hertz problem with large deformations} \label{subsec:large deformation}

Finally, in this section the two-dimensional frictionless Hertz problem is studied considering large deformations and strains.
For that purpose, a Neo-Hookean material constitutive law {(an hyper-elastic law that considers finite strains)} with Young's modulus $E = 1$ and Poisson's ratio $\nu = 0.3$,
has been used for the deformable body.

\begin{figure}[!ht]
  \centering
  \subfigure[Stress magnitude distribution for the undeformed mesh with $N_2 / S_0$ method.]{\label{fig:mesh4_b} \includegraphics[width= 8.5cm]{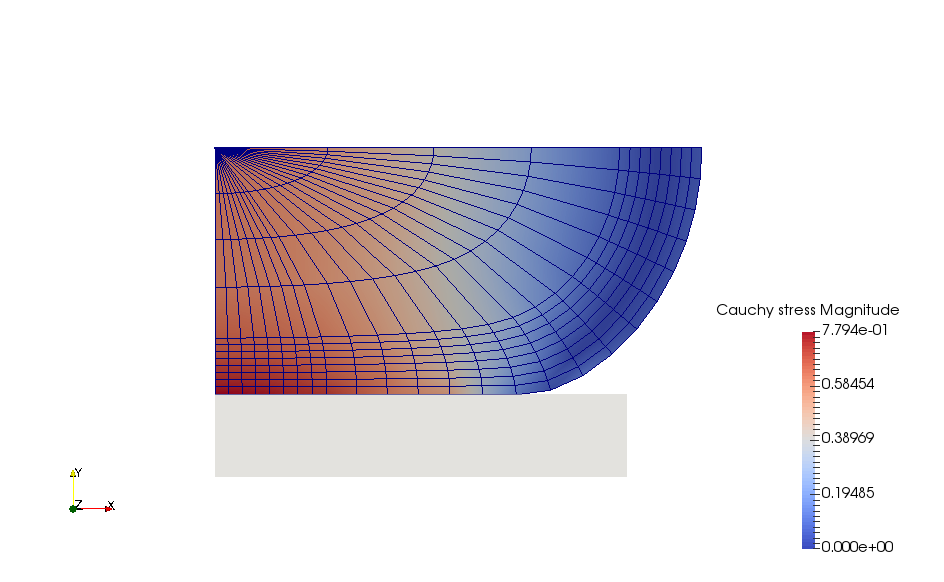}}\hfill
  \caption{2D large deformation Hertz contact problem with a uniform downward displacement $u_y = -0.4$.}
  \label{fig:mesh4}
\end{figure}

\begin{figure}[!ht]
  \centering
  \subfigure[Reference numerical contact pressure with $N_2 / S_0$ and $N_3 / S_1$ methods.]{\label{fig:mesh4_b mult ref} \includegraphics[width= 7cm]{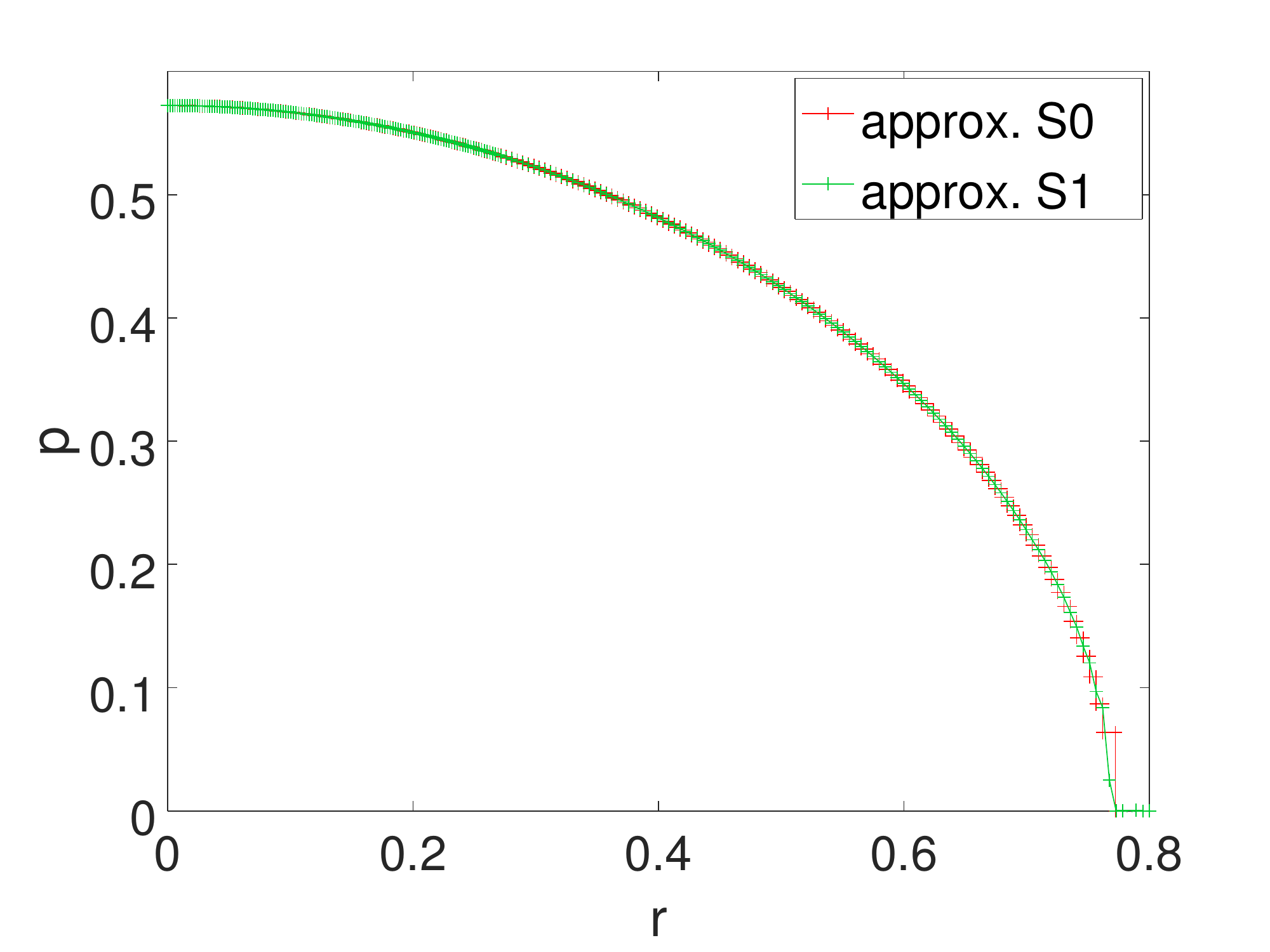}}
  \subfigure[Numerical contact pressure with $N_2 / S_0$ and $N_3 / S_1$ methods.]{\label{fig:mesh4_b mult} \includegraphics[width= 7cm]{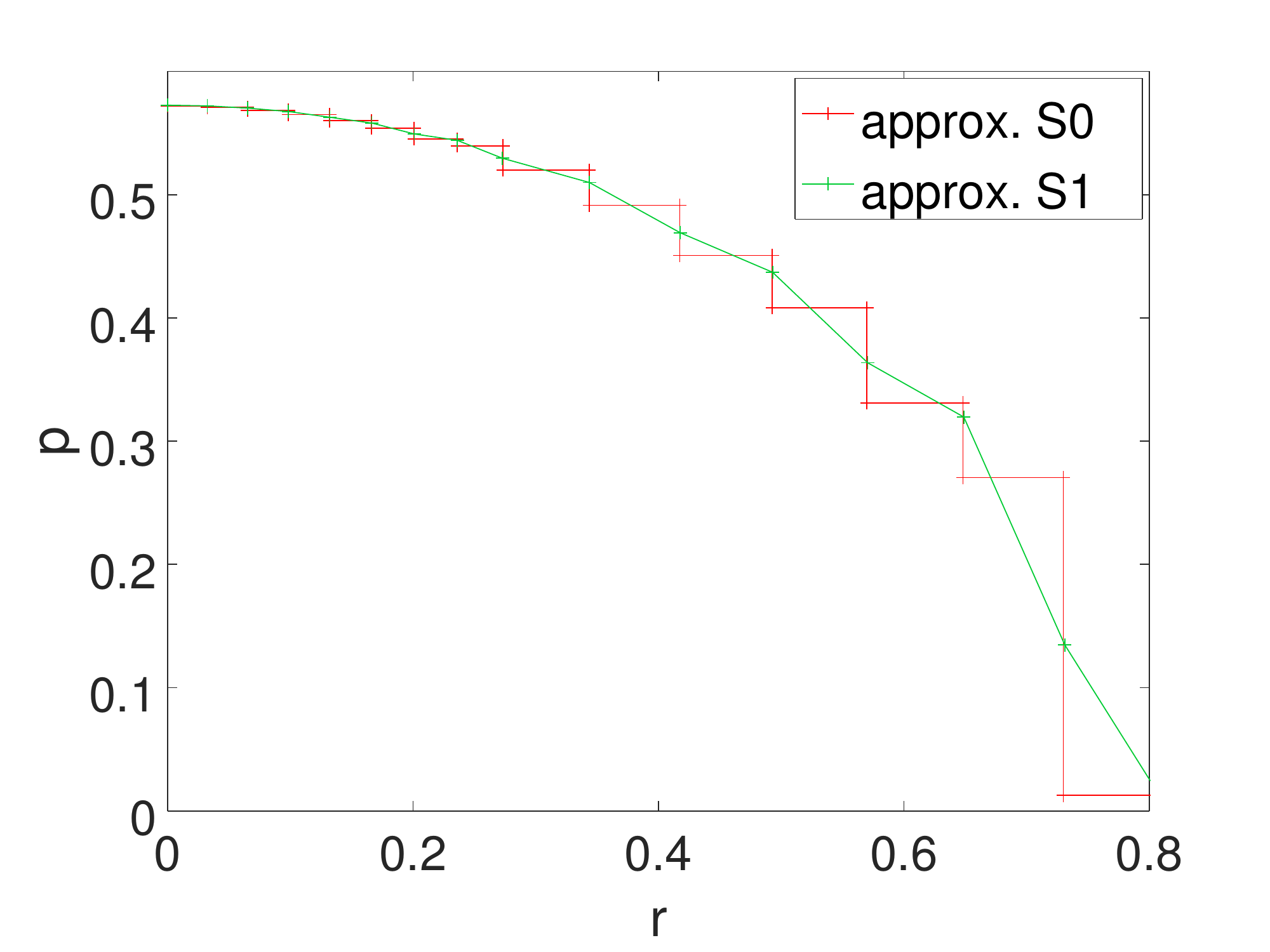}}
  \caption{2D large deformation Hertz contact problem with a uniform downward displacement $u_y = -0.4$.}
  \label{fig:mesh4_mult}
\end{figure}

As in Section \ref{subsec:Hertz problems 2d}, the performance
of the $N_2 / S_0$ and $N_3 / S_1$ method are analysed and the problem is modelled as an elastic quarter of disc but modifying its boundary conditions:
instead of pressure, a uniform downward displacement $u_y=-0.4$ is applied on the top surface (see Figure \ref{fig:mesh4}). In this large deformation framework the exact solution is unknown: the error of the computed displacement and Lagrange multipliers are studied
taking a refined numerical solution as reference. The large deformation of the body and computed contact pressure are presented in Figure \ref{fig:mesh4_mult}.

As in the previous case optimal results are obtained for
the computed displacement and Lagrange multipliers (see Figure \ref{fig:disp 0point4}).
\begin{figure}[!ht]
  \centering
   \subfigure[$N_2 / S_0$ method.]{\label{fig:sec3 mult dir 0.4} 
\begin{tikzpicture}
\begin{loglogaxis}[width=7.3cm ,ymin=3e-07,xmin=1.1e-02,xmax=2.7e0,
xlabel={h},ylabel={absolute errors},legend pos={south east},grid=major]

\addplot table[x=h,y=L2_abs] {figure_latex_figure_data_large_dirichlet_data_large_dirichlet.res};
\addplot table[x=h,y=H1_abs] {figure_latex_figure_data_large_dirichlet_data_large_dirichlet.res};
\addplot[brown,mark=pentagon*] table[x=h,y=L2_mult_abs] {figure_latex_figure_data_large_dirichlet_data_large_dirichlet.res};

\slopeTriangleAbove{0.03}{0.045}{1.4e-06}{3.711e-06}{2.404}{blue}; 
\slopeTriangleAbove{0.03}{0.045}{0.0002}{0.000364}{1.476}{red}; 
\slopeTriangleBelow{0.03}{0.045}{0.007}{0.0101}{0.908}{brown}; 

\legend{ \footnotesize $L^2$-norm, \footnotesize $H^1$-{semi-}norm,\footnotesize mult. $L^2$-norm}
\end{loglogaxis}
\end{tikzpicture}
} \hfill
     \subfigure[$N_3 / S_1$ method.]{\label{fig:sec3 disp dir 0.4} 
 \begin{tikzpicture}
\begin{loglogaxis}[width=7.3cm ,ymin=2e-07,xmin=2e-02,xmax=2.7e0,
xlabel={h},ylabel={absolute errors},legend pos={south east},grid=major]

\addplot table[x=h,y=L2_abs] {figure_latex_figure_data_large_dirichlet_data_large_dirichlet_P3.res};
\addplot table[x=h,y=H1_abs] {figure_latex_figure_data_large_dirichlet_data_large_dirichlet_P3.res};
\addplot[brown,mark=pentagon*] table[x=h,y=L2_mult_abs] {figure_latex_figure_data_large_dirichlet_data_large_dirichlet_P3.res};

\slopeTriangleAbove{0.055}{0.08}{0.7e-06}{2.1990e-06}{3.055}{blue}; 
\slopeTriangleAbove{0.055}{0.08}{0.0002}{4.2061e-04}{1.984}{red}; 
\slopeTriangleBelow{0.055}{0.08}{0.002}{0.0035990}{1.568}{brown}; 

\legend{ \footnotesize $L^2$-norm, \footnotesize $H^1$-{semi-}norm,\footnotesize mult. $L^2$-norm}
\end{loglogaxis}
\end{tikzpicture}           
}
  \caption{2D large deformation Hertz contact problem method with a uniform downward displacement $u_y = -0.4$.
   Absolute displacement errors in $L^2$-norm and $H^1$-{semi-}norm and Lagrange multipliers error in $L^2$-norm.} 
  \label{fig:disp 0point4}
\end{figure}
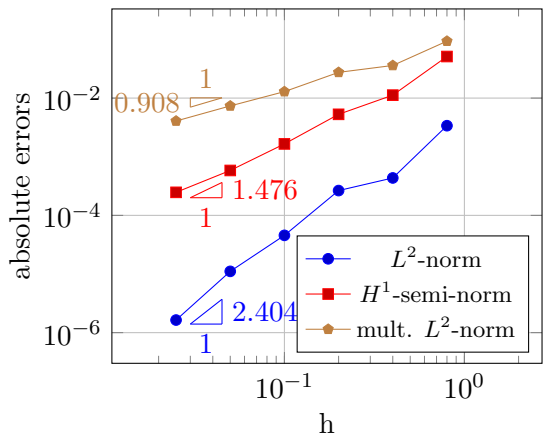
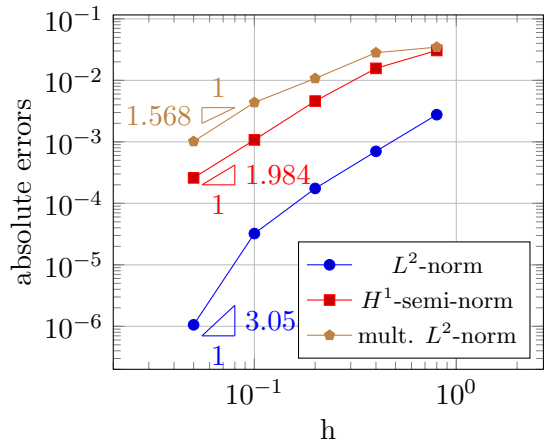\\

\section*{Conclusions}
\label{sec: conclusion}
\addcontentsline{toc}{section}{Conclusions}

In this work, we presented an {optimal} \textit{a priori} error estimate of frictionless unilateral contact problems between a deformable body and a rigid one for an augmented lagrangian method. 

For the numerical point of view, we observe an {optimality} of this method for both variables, the displacement and the Lagrange multipliers. In our experiments, we used a NURBS of degree $2$ for the primal space and B-Spline of degree $0$ for the dual space as well as a NURBS of degree $3$ for the primal space and B-Spline of degree $1$ for the dual space. Thanks to this choice of approximation spaces, we observe a stability of the Lagrange multipliers, indeed no oscillation are observed, and a well approximation of the pressure in the two- and tree-dimensional case. 

\section*{Acknowledgements}
\label{sec: acknowledgements}
\addcontentsline{toc}{section}{acknowledgements}
This work has been partially supported by Michelin under the contract A10-4087.

\newpage

\addcontentsline{toc}{section}{Appendices}
\section*{{Appendix 1.}}
\label{sec:appendix1}

In this appendix, we provide the ingredients needed to fully discretise the problem \eqref{eq:augmultdiscrete} as well as its large deformation version that we have used in Section \ref{sec:sec3}. First we introduce the contact status, an active-set strategy for the discrete problem, and then the fully discrete problem. For the purpose of this appendix, we take notations suitable to large deformation and denote by $g_n$ the distance between the rigid and the deformable body. In small deformation, it holds $g_n(u)=u \cdot n$.
\subsection*{Contact status}
\label{subsec:contact stat}
Let us first deal with the contact status. The active-set strategy is defined in \cite{hueber-Wohl-05,hueber-Stadler-Wohlmuth-08} and is updated at each iteration of Newton. Due to the deformation, parts of the workpiece may come into contact or conversely may loose contact. This change of contact status changes the loading that is applied on the boundary of the mesh. This method is used to track the location of contact during the change in boundary conditions.\\
An equivalent formulation of the contact status presented in \cite{fabre-a-priori-iga-17}, let $K$ be a control point of the B-Spline space \eqref{def:space_mult}, let $(\Pi^h_\lambda \cdot )_K$ be the local projection defined in \eqref{def:l2proj_K} and in the same way as previously done:
\begin{itemize}
\item if $\lambda_{n,K} {B}_{K} + r (\Pi_\lambda g_n)_K <0$, the control point $K$ is active;
\item if $\lambda_{n,K} {B}_{K} + r (\Pi_\lambda g_n)_K \geq 0$, the control point $K$ is inactive.
\end{itemize}

\subsection{Discrete Problem}
\label{subsec: discrete problem alm }

Regarding the augmented Lagrange multiplier method, the contact contribution of the work is expressed as follows 
\begin{eqnarray}
\begin{array}{l}
\dis  W_{c,r} = \frac{1}{2r}  \int_{\Gamma_C} [ \lambda + r g_n]_-^2 - \lambda^2 \ \textrm{d}\Gamma .
\end{array}\nonumber
\end{eqnarray}
The contact contribution of the virtual work is expressed as follows 
\begin{eqnarray}
\begin{array}{l}
\dis  \delta W_{c,r} = - \int_{\Gamma_C} [ \lambda + r g_n]_- \delta g_n \ \textrm{d}\Gamma - \frac{1}{r}  \int_{\Gamma_C} (\lambda +[ \lambda + r g_n]_- ) \delta \lambda \ \textrm{d}\Gamma .
\end{array}\nonumber
\end{eqnarray}

\noindent In order to implement this method, we need to use our local gap. For simplification, the discretised contact contribution can be expressed as follows
\begin{eqnarray}
\begin{array}{ll}
\dis  \delta W_{c,r} &\dis = - \int_{\Gamma_C} \left[  \sum_{K=1}^{n_{\lambda}} \lambda_{n,K} {B}_{K} + r (\Pi_\lambda g_n)_K {B}_{K} \right]_- \sum_{A=1}^{n_{u}} \delta u_A N_A  n \ \textrm{d}\Gamma \\
&\dis - \frac{1}{r}  \int_{\Gamma_C} \left( \sum_{K=1}^{n_{\lambda}} \lambda_{n,K} {B}_{K}  + \left[  \sum_{K=1}^{n_{\lambda}} \lambda_{n,K} {B}_{K} + r (\Pi_\lambda g_n)_K {B}_{K} \right]_- \right) \sum_{\kappa=1}^{n_{\lambda}} \delta \lambda_{n,\kappa} {B}_{\kappa} \ \textrm{d}\Gamma.
\end{array}\nonumber
\end{eqnarray}
\noindent Now, we can distinguish between the active part and the inactive part, it holds:
\begin{eqnarray}
\begin{array}{ll}
\dis  \delta W_{c,r} &\dis = \int_{\Gamma_C} \Big{(}  \sum_{K, act} \lambda_{n,K} {B}_{K} + r (\Pi_\lambda g_n)_K  {B}_{K}\Big{)}  \sum_{A} \delta u_A N_A  n \ \textrm{d}\Gamma \\
&\dis - \frac{1}{r}  \int_{\Gamma_C} \left( \sum_{K} \lambda_{n,K} {B}_{K}  -  \Big{(}  \sum_{K, act} \lambda_{n,K} {B}_{K} + r (\Pi_\lambda g_n)_K {B}_{K} \Big{)} \right) \sum_{\kappa} \delta \lambda_{n,\kappa} {B}_{\kappa}\ \textrm{d}\Gamma ,\\[0.4 cm]
&\dis =  \sum_{K, act} \sum_{A}  \int_{\Gamma_C} \Big{(}  \lambda_{n,K} {B}_{K} + r (\Pi_\lambda g_n)_K {B}_{K} \Big{)}   \delta u_A N_A  n \ \textrm{d}\Gamma \\
&\dis - \frac{1}{r} \sum_{K,inact}  \sum_{\kappa} \int_{\Gamma_C}  \lambda_{n,K} {B}_{K} \delta \lambda_{n,\kappa} {B}_{\kappa}\ \textrm{d}\Gamma  + \sum_{K, act} \sum_{\kappa} \int_{\Gamma_C}  (\Pi_\lambda g_n)_K {B}_{K}    \delta \lambda_{n,\kappa} {B}_{\kappa} \ \textrm{d}\Gamma , \\[0.4 cm]
&\dis =  \sum_{K, act} \sum_{A} \delta u_A \int_{\Gamma_C} \Big{(}  \lambda_{n,K} {B}_{K} + r (\Pi_\lambda g_n)_K {B}_{K} \Big{)}    N_A  n \ \textrm{d}\Gamma \\
&\dis - \frac{1}{r} \sum_{K,inact}  \sum_{\kappa} \delta \lambda_{n,\kappa} \int_{\Gamma_C}  \lambda_{n,K} {B}_{K}  {B}_{\kappa}\ \textrm{d}\Gamma  + \sum_{K, act} \sum_{\kappa}  \delta \lambda_{n,\kappa} \int_{\Gamma_C}  (\Pi_\lambda g_n)_K {B}_{K} {B}_{\kappa} \ \textrm{d}\Gamma , \\[0.4 cm]
&\dis =  \sum_{K, act}  \delta \bfu^T \int_{\Gamma_C} \Big{(}  \lambda_{n,K} {B}_{K} + r (\Pi_\lambda g_n)_K {B}_{K} \Big{)}    \bfN \ \textrm{d}\Gamma \\
&\dis - \frac{1}{r} \sum_{K,inact}  \delta \bflambda_{n}^T \int_{\Gamma_C}  \lambda_{n,K} {B}_{K}  {\bfB_\lambda}\ \textrm{d}\Gamma  + \sum_{K, act}  \delta \bflambda_{n}^T \int_{\Gamma_C}  (\Pi_\lambda g_n)_K {B}_{K} {\bfB_\lambda} \ \textrm{d}\Gamma .
\end{array}\nonumber
\end{eqnarray}
The residual for Newton-Raphson iterative scheme is obtained as $$\dis R_r = \begin{bmatrix} R_u + R_{u,r} \\ R_\lambda + R_{\lambda,r} \end{bmatrix}= \begin{bmatrix} \dis  \int_{\Gamma_C} \sum_{K, act} ( \lambda_{n,K} {B}_{K} ) \bfN  \ \textrm{d}\Gamma    + r \int_{\Gamma_C} \sum_{K, act} ( (\Pi_\lambda g_n)_K {B}_{K}   ) \bfN \ \textrm{d}\Gamma \\ \dis \bfN_{\lambda,g}  -\frac{1}{r} \int_{\Gamma_C}( \sum_{K,inact} \lambda_{n,K} {B}_{K} ) B_\lambda \ \textrm{d}\Gamma \end{bmatrix} .$$
\noindent The linearization and active set strategy yield 
\begin{eqnarray}
\begin{array}{ll}
\dis \Delta \delta W_{c,r} &\dis =   \delta \bfu^T \int_{\Gamma_C} \sum_{K, act}  \Big{(} \Delta \lambda_{n,K} {B}_{K} + r (\Pi_\lambda \Delta g_n)_K {B}_{K} \Big{)}    \bfN \ \textrm{d}\Gamma \\
&\dis - \frac{1}{r}  \delta \bflambda_{n}^T \int_{\Gamma_C}\sum_{K,inact}   \Delta  \lambda_{n,K} {B}_{K}  {\bfB_\lambda}\ \textrm{d}\Gamma  +   \delta \bflambda_{n}^T \int_{\Gamma_C} \sum_{K, act}  (\Pi_\lambda \Delta g_n)_K {B}_{K} {\bfB_\lambda} \ \textrm{d}\Gamma .
\end{array}\nonumber
\end{eqnarray}
We define the following matrix $$\bfB_{\lambda,act} = \begin{bmatrix}  B_{1,act}(\zeta) \\   \vdots \\  B_{n_{\lambda},act}(\zeta)  \end{bmatrix}, \quad \bfB_{\lambda,inact} = \begin{bmatrix}  B_{1,inact}(\zeta) \\   \vdots \\  B_{n_{\lambda},inact}(\zeta)  \end{bmatrix}.$$
The discretised of contact contribution can be expressed as follows
\begin{eqnarray}
\begin{array}{ll}
\dis \Delta \delta W_{c,r} &\dis =   \delta \bfu^T \int_{\Gamma_C} \bfN \bfB_{\lambda,act}^T \ \textrm{d}\Gamma \Delta \bflambda_{n} +  \delta \bfu^T \Big{(}\sum_{K, act}  \frac{r}{\int_{\Gamma_C} {B}_{K} \ \textrm{d}\Gamma}  \int_{\Gamma_C} \bfN {B}_{K} \ \textrm{d}\Gamma \int_{\Gamma_C} {B}_{K} \bfN^T   \ \textrm{d}\Gamma \Big{)} \Delta \bfu \\
&\dis  +   \delta \bflambda_{n}^T \Big{(}\sum_{K, act}  \frac{1}{\int_{\Gamma_C} {B}_{K} \ \textrm{d}\Gamma}  \int_{\Gamma_C} \bfB {B}_{K} \ \textrm{d}\Gamma \int_{\Gamma_C} {B}_{K} \bfN^T   \ \textrm{d}\Gamma \Big{)} \Delta \bfu - \frac{1}{r}  \delta \bflambda_{n}^T \int_{\Gamma_C} \bfB_{\lambda} \bfB_{\lambda,inact}^T \ \textrm{d}\Gamma  \Delta  \bflambda_{n,K} .
\end{array}\nonumber
\end{eqnarray}

 \addcontentsline{toc}{section}{Bibliography}
  
\bibliographystyle{siam}

\begin{thebibliography}{10}

\bibitem{al-cu1988}
{\sc P.~Alart and A.~Curnier}, {\em A generalized {N}ewton method for contact
  problems with friction}, Journal de M\'ecanique Th\'eorique et Appliqu\'ee, 7 (1988), pp.~67--82.  
  
    
\bibitem{fabre-a-priori-iga-17}  
{\sc P.~Antolìn, A.~Buffa and M.~Fabre}, {\em A priori error for unilateral contact problems with Lagrange multiplier and IsoGeometric Analysis},
submitted in IMA Journal of Numerical Analysis, (2017).
     
\bibitem{daveiga06}
{\sc Y.~Bazilevs, L.~Beir\~ao~da Veiga, J.~A. Cottrell, T.~J.~R. Hughes, and
  G.~Sangalli}, {\em Isogeometric analysis: Approximation, stability and error
  estimates for h-refined meshes}, Mathematical Models and Methods in Applied
  Sciences, 16 (2006), pp.~1031--1090.
  
\bibitem{buffa09}
{\sc L.~Beir\~ao~da Veiga, A.~Buffa, J.~Rivas and G.~Sangalli}, {\em Some estimates for $h-p-k$-refinement in Isogeometric Analysis}, Numerische Mathematik, 118 (2011), pp.~271--305.  

\bibitem{daveiga2014}
{\sc L.~Beir\~ao~da~Veiga, A.~Buffa, G.~Sangalli, and R.~V\'{a}zquez}, {\em
  Mathematical analysis of variational isogeometric methods}, Acta Numerica, 23
  (2014), pp.~157--287.

\bibitem{ben-belgacem-renard-03}
{\sc F.~Ben~Belgacem and Y.~Renard}, {\em Hybrid finite element methods for the
  {S}ignorini problem}, Mathematics of Computation, 72 (2003), pp.~1117--1145.
  
\bibitem{bica-kozia-08}
{\sc T.~K. Bi\'cani\'c}, {\em Semismooth newton method for frictional contact
  between pseudo-rigid bodies}, Computer Methods in Applied Mechanics and
  Engineering, 197 (2008).  
  
  
  \bibitem{burman-16}
{\sc E.~Burman, P.~Hansbo and M.~Larson}, {\em Augmented Lagrangian finite element methods for contact problems},
Submitted (2016).

\bibitem{fortinbook13}
{\sc D.~Boffi, F.~Brezzi, and M.~Fortin}, {\em Mixed finite element methods and
  applications}, Computational Mathematics, Springer, 2013.

\bibitem{brivadis15}
{\sc E.~Brivadis, A.~Buffa, B.~Wohlmuth, and L.~Wunderlich}, {\em Isogeometric
  mortar methods}, Computer Methods in Applied Mechanics and Engineering, 284
  (2015), pp.~292--319.

\bibitem{coorevits-hild-lhalouani-sassi-02}
{\sc P.~Coorevits, P.~Hild, K.~Lhalouani and T.~Sassi, }, {\em Mixed finite element methods for unilateral problems: convergence analysis and numerical studies}, Mathematics of Computation, 71 (2002), pp.~67--82.    

\bibitem{franz-overview}
{\sc F.~Chouly, M.~Fabre, P.~Hild, R.~Mlika, J.~Pousin and Y.~Renard}, {\em An Overview of Recent Results on Nitsche's Method for Contact Problems},
Geometrically Unfitted Finite Element Methods and Application, Springer International Publishing, 2017.

\bibitem{chouly-hild-2013}
{\sc F.~Chouly and P.~Hild}, {\em A Nitsche-Based Method for Unilateral Contact Problems: Numerical Analysis},
SIAM Journal on Numerical Analysis, 51 (2013), pp.~1295--1307.

\bibitem{chouly-hild-renard-2014}
{\sc F.~Chouly, P.~Hild and Y.~Renard}, {\em Symmetric and non-symmetric variants of Nitsche's method for contact problems in elasticity: Theory and numerical experiments},
Mathematics of Computation, 84 (2015), pp.~1089--1112.

\bibitem{chouly-hild-renard-2015-1}
{\sc F.~Chouly, P.~Hild and Y.~Renard}, {\em A Nitsche finite element method for dynamic contact : 1. Semi-discrete problem analysis and time-marching schemes},
ESAIM: Mathematical Modelling and Numerical Analysis, 49 (2015), pp.~481--502.

\bibitem{chouly-hild-renard-2015-2}
{\sc F.~Chouly, P.~Hild and Y.~Renard}, {\em A Nitsche finite element method for dynamic contact : 2. Semi-discrete problem analysis and time-marching schemes},
ESAIM: Mathematical Modelling and Numerical Analysis, 49 (2015), pp.~503--528.
  
\bibitem{lorenzis15}
{\sc L.~De~Lorenzis, J.~Evans, T.~Hughes, and A.~Reali}, {\em Isogeometric
  collocation: Neumann boundary conditions and contact}, Computer Methods in
  Applied Mechanics and Engineering, 284 (2015), pp.~21--54.

\bibitem{lorenzisreview}
{\sc L.~De~Lorenzis, P.~Wriggers, and T.~J. Hughes}, {\em Isogeometric contact:
  A review}, GAMM-Mitteilungen, 37 (2014), pp.~85--123.

\bibitem{lorenzis11}
{\sc L.~De~Lorenzis, P.~Wriggers, and G.~Zavarise}, {\em A large deformation frictional contact formulation using NURBS-based isogeometric analysis}, International Journal for Numerical Methods in Engineering (2011).

\bibitem{lorenzis12}
{\sc L.~De~Lorenzis, P.~Wriggers, and G.~Zavarise}, {\em A mortar formulation
  for 3d large deformation contact using nurbs-based isogeometric analysis and
  the augmented lagrangian method}, Springer-Verlag, 49 (2012), pp.~1--20.

\bibitem{drouet15}
{\sc G.~Drouet and P.~Hild}, {\em Optimal convergence for discrete variational
  inequalities modelling signorini contact in 2d and 3d without additional
  assumptions on the unknown contact set}, SIAM Journal on Numerical Analysis, 53
  (2015), pp.~1488--1507.
  
  \bibitem{fabre-pousin-renard-2016}
{\sc M.~Fabre, J.~Pousin and Y.~Renard}, {\em A fictitious domain method for frictionless contact problems in elasticity using Nitsche's method},
SMAI Journal of Computational Mathematics, 2 (2016), pp.19--50.

\bibitem{haslinger-96}
{\sc J.~Haslinger, I.~Hlav\'a{\v{c}}ek, and J.~Ne{\v{c}}as}, {\em Handbook of
  Numerical Analysis (eds. P.G. Ciarlet and J.L. Lions)}, vol.~IV, North
  Holland, 1996, ch.~2. ``Numerical methods for unilateral problems in solid
  mechanics'', pp.~313--385.

 \bibitem{Hertz1882}
{\sc H.~Hertz}, {\em \"Ueber die ber\"uhrung fester elastischer k\"orper}, Journal f\"ur die reine und angewandte Mathematik,  92 (1882), pp.~156--171.

\bibitem{hild-laborde-02}
{\sc P.~Hild and P.~Laborde}, {\em Quadratic finite element methods for
  unilateral contact problems}, Applied Numerical Mathematics, 41 (2002), pp.~401--421.

\bibitem{hueber-Stadler-Wohlmuth-08}
{\sc S.~H{\"u}eber, G.~Stadler, and B.~I. Wohlmuth}, {\em A primal-dual active
  set algorithm for three-dimensional contact problems with coulomb friction},
  SIAM Journal on Scientific Computing, 30 (2008), pp.~572--596.

\bibitem{hueber-Wohl-05}
{\sc S.~H{\"u}eber and B.~I. Wohlmuth}, {\em A primal--dual active set strategy
  for non-linear multibody contact problems}, Computer Methods in Applied
  Mechanics and Engineering, 194 (2005).

\bibitem{hughes05}
{\sc T.~J.~R. Hughes, J.~A. Cottrell, and Y.~Bazilev}, {\em Isogeometric
  analysis: Cad, finite elements, nurbs, exact geometry and mesh refinement},
  Computer Methods in Applied Mechanics and Engineering, 194 (2005), pp.~4135--4195.

\bibitem{johnsonbook}
{\sc K.~L. Johnson}, {\em Contact mechanics}, Cambridge University Press, 1985.

\bibitem{kon-sch-13}
{\sc A.~Konyukhov and K.~Schweizerhof}, {\em Computational contact mechanics},
  vol.~67, Applied and Computational Mechanics, 2013.

\bibitem{laursen-02}
{\sc T.~A. Laursen}, {\em Computational contact and impact mechanics},
  Springer-Verlag, Berlin, 2003.

\bibitem{lions-magenes-72}
{\sc J.~L. Lions and E.~Magenes}, {\em Non-homogeneous boundary value problems
  and applications}, Springer-Verlag, Berlin, New York, 1972.
  
  \bibitem{chouly-mlika-renard-2017}
{\sc R.~Mlika, Y.~Renard and F.~Chouly}, {\em An unbiased Nitsche's formulation of large deformation frictional contact and self-contact},
Computer Methods in Applied Mechanics and Engineering, 325 (2017), pp.~265--288.


\bibitem{Moussaoui92}
{\sc M.~Moussaoui and K.~Khodja}, {\em R{\'e}gularit{\'e} des solutions d'un
  probl{\`e}me m{\^e}l{\'e} dirichlet-signorini dans un domaine polygonal
  plan}, Communications in Partial Differential Equations,  (1992), pp.~805--826.

\bibitem{Pauletti2015}
{\sc M.~S. Pauletti, M.~Martinelli, N.~Cavallini, and P.~Antolin}, {\em
  {Igatools: An isogeometric analysis library}}, SIAM Journal on Scientific Computing, 37
  (2015).

\bibitem{wohl12}
{\sc A.~Popp, B.~I. Wohlmuth, M.~W. Gee, and W.~A. Wall}, {\em Dual quadratic
  mortar finite element methods for 3d finite deformation contact}, SIAM
  Journal on Scientific Computing, 34 (2012), pp.~B421--B446.

\bibitem{renard-poulios14_2}
{\sc K.~Poulios and Y.~Renard}, {\em An unconstrained integral approximation of
  large sliding frictional contact between deformable solids}, Computers and
  Structures, 153 (2015), pp.~75--90.
  
  \bibitem{renard-2012b}
{\sc Y.~Renard}, {\em {Generalized Newton's methods for the approximation and
  resolution of frictional contact problems in elasticity}}, Comp. Methods
  Appl. Mech. Engrg., 256 (2013), pp.~38--55.

\bibitem{seitz16}
{\sc A.~Seitz, P.~Farah, J.~Kremheller, B.~I. Wohlmuth, W.~A. Wall, and
  A.~Popp}, {\em Isogeometric dual mortar methods for computational contact
  mechanics}, Computer Methods in Applied Mechanics and Engineering, 301
  (2016), pp.~259--280.

\bibitem{Temizer11}
{\sc I.~Temizer, P.~Wriggers, and T.~Hughes}, {\em Contact treatment in
  isogeometric analysis with NURBS}, Computer Methods in Applied Mechanics
  and Engineering, 200 (2011), pp.~1100--1112.

\bibitem{Temizer12}
{\sc I.~Temizer, P.~Wriggers, and T.~Hughes}, {\em Three-dimensional
  mortar-based frictional contact treatment in isogeometric analysis with
  NURBS}, Computer Methods in Applied Mechanics and Engineering, 209--212
  (2012), pp.~115--128.

\bibitem{Wohlmuth00}
{\sc B.~I. Wohlmuth}, {\em A mortar finite element method using dual spaces for
  the lagrange multiplier}, SIAM Journal on Numerical Analysis, 38 (2000),
  pp.~989--1012.

\bibitem{wriggers-06}
{\sc P.~Wriggers}, {\em Computational contact mechanics (Second Edition)},
  Wiley, 2006.

\bibitem{lorenzis9}
{\sc G.~Zavarise and L.~D. Lorenzis}, {\em The node-to-segment algorithm for 2d
  frictionless contact: Classical formulation and special cases}, Computer Methods in Applied Mechanics and Engineering,  (2009).

\end{thebibliography}


\end{document}